\documentclass[10pt,a4paper]{article}
\usepackage[utf8]{inputenc}
\usepackage[T1]{fontenc}
\usepackage{geometry}
\usepackage{amsthm}
\usepackage[francais,english]{babel}
\usepackage{amssymb}
\usepackage{lmodern}
\usepackage{amsmath,amsfonts}
\usepackage{mathrsfs} 
\usepackage{dsfont}
\usepackage{stmaryrd}
 
\usepackage{graphicx}
\usepackage{fullpage}
\usepackage[nottoc,notlot,notlof]{tocbibind}
\usepackage[colorlinks=true,urlcolor=black,linkcolor=black,citecolor=black]{hyperref}
\numberwithin{equation}{section}
\usepackage{mathtools,array}
\usepackage{verbatimbox}
\usepackage{color}
\usepackage{bm}
\usepackage{tikz-cd}
\usepackage{pst-node}
\usetikzlibrary{matrix}
\usepackage{mathtools}
\usepackage{verbatimbox}
\usepackage{diagbox}
\usepackage{array}
\newcolumntype{C}{>{$\displaystyle} c <{$}}
\usepackage{makecell}
\usepackage{float}
\usepackage{tabu}
\usepackage{cancel}
\usepackage{pdflscape}
\makeatletter
\def\env@dmatrix{\hskip -\arraycolsep
	\let\@ifnextchar\new@ifnextchar
	\def\arraystretch{2}%
	\array{*{\c@MaxMatrixCols}{>{\displaystyle}c}}%
}

\makeatother

\begin{document}

	\renewcommand{\thefootnote}{\fnsymbol{footnote}}
	
	\title{
On the Moduli Space of Null Curves in Klein's Quadric}
	\author{Alexis Michelat\footnote{Department of Mathematics, ETH Zentrum, CH-8093 Zürich, Switzerland.}}
	\date{\today}
	
	\maketitle
	
	\vspace{-0.5em}

	\begin{abstract}
		We study the moduli space of null curves in Klein's quartic in the four-dimensional (complex) projective plane using methods developed by Robert Bryant. As a consequence, we show that minimal surfaces with $9$ embedded planar ends do not exist and formulate some conjectures about the previous moduli space.
	\end{abstract}
	
	\tableofcontents
	\vspace{0cm}
	\begin{center}
	\end{center}

	\theoremstyle{plain}
	\newtheorem*{theorem*}{Theorem}
	\newtheorem{theorem}{Theorem}[section]
	\newenvironment{theorembis}[1]
	{\renewcommand{\thetheorem}{\ref{#1}$'$}%
		\addtocounter{theorem}{-1}%
		\begin{theorem}}
		{\end{theorem}}
	\renewcommand*{\thetheorem}{\Alph{theorem}}
	\newtheorem{lemme}[theorem]{Lemma}
	\newtheorem{propdef}[theorem]{Définition-Proposition}
	\newtheorem{prop}[theorem]{Proposition}
	\newtheorem{cor}[theorem]{Corollary}
	\theoremstyle{definition}
	\newtheorem*{definition}{Definition}
	\newtheorem{defi}[theorem]{Definition}
	\newtheorem{rem}[theorem]{Remark}
	\newtheorem*{remast}{Remark}
	\newtheorem*{conjecture}{Conjecture}
	\newtheorem{rems}[theorem]{Remarks}
	\newtheorem{exemple}[theorem]{Example}
	\renewcommand\hat[1]{%
		\savestack{\tmpbox}{\stretchto{%
				\scaleto{%
					\scalerel*[\widthof{\ensuremath{#1}}]{\kern-.6pt\bigwedge\kern-.6pt}%
					{\rule[-\textheight/2]{1ex}{\textheight}}
				}{\textheight}%
			}{0.5ex}}%
		\stackon[1pt]{#1}{\tmpbox}
	}
	\parskip 1ex
	\newcommand{\totimes}{\ensuremath{\,\dot{\otimes}\,}}
	\newcommand{\vc}[3]{\overset{#2}{\underset{#3}{#1}}}
	\newcommand{\conv}[1]{\ensuremath{\underset{#1}{\longrightarrow}}}
	\newcommand{\A}{\ensuremath{\vec{A}}}
	\newcommand{\B}{\ensuremath{\vec{B}}}
	\newcommand{\C}{\ensuremath{\mathbb{C}}}
	\newcommand{\D}{\ensuremath{\nabla}}
	\newcommand{\E}{\ensuremath{\vec{E}}}
	\newcommand{\Q}{\ensuremath{\vec{Q}}}
	\newcommand{\z}{\ensuremath{\bar{z}}}
	\newcommand{\R}{\ensuremath{\mathbb{R}}}
	\renewcommand{\P}{\ensuremath{\mathbb{P}}}
	\newcommand{\N}{\ensuremath{\mathbb{N}}}
	\newcommand{\Z}{\ensuremath{\mathbb{Z}}}
	\newcommand{\p}[1]{\ensuremath{\partial_{#1}}}
	\newcommand{\Res}{\ensuremath{\mathrm{Res}}}
	\newcommand{\lp}[2]{\ensuremath{\mathrm{L}^{#1}(#2)}}
	\renewcommand{\wp}[3]{\ensuremath{\left\Vert #1\right\Vert_{\mathrm{W}^{#2}(#3)}}}
	\newcommand{\np}[3]{\ensuremath{\left\Vert #1\right\Vert_{\mathrm{L}^{#2}(#3)}}}
	\newcommand{\nc}[3]{\ensuremath{\left\Vert #1\right\Vert_{C^{#2}(#3)}}}
	\newcommand{\h}{\ensuremath{\vec{h}}}
	\renewcommand{\Re}{\ensuremath{\mathrm{Re}\,}}
	\renewcommand{\Im}{\ensuremath{\mathrm{Im}\,}}
	\newcommand{\diam}{\ensuremath{\mathrm{diam}\,}}
	\newcommand{\leb}{\ensuremath{\mathscr{L}}}
	\newcommand{\supp}{\ensuremath{\mathrm{supp}\,}}
	\renewcommand{\phi}{\ensuremath{\vec{\Phi}}}
	\renewcommand{\H}{\ensuremath{\vec{H}}}
	\renewcommand{\L}{\ensuremath{\mathrm{L}}}
	\renewcommand{\epsilon}{\ensuremath{\varepsilon}}
	\renewcommand{\bar}{\ensuremath{\overline}}
	\newcommand{\s}[2]{\ensuremath{\langle #1,#2\rangle}}
	\newcommand{\pwedge}[2]{\ensuremath{\,#1\wedge#2\,}}
	\newcommand{\bs}[2]{\ensuremath{\left\langle #1,#2\right\rangle}}
	\newcommand{\scal}[2]{\ensuremath{\langle #1,#2\rangle}}
	\newcommand{\sg}[2]{\ensuremath{\left\langle #1,#2\right\rangle_{\mkern-3mu g}}}
	\newcommand{\n}{\ensuremath{\vec{n}}}
	\newcommand{\ens}[1]{\ensuremath{\left\{ #1\right\}}}
	\newcommand{\lie}[2]{\ensuremath{\left[#1,#2\right]}}
	\newcommand{\g}{\ensuremath{g}}
	\newcommand{\e}{\ensuremath{\vec{e}}}
	\newcommand{\ig}{\ensuremath{|\vec{\mathbb{I}}_{\phi}|}}
	\newcommand{\ik}{\ensuremath{\left|\mathbb{I}_{\phi_k}\right|}}
	\newcommand{\w}{\ensuremath{\vec{w}}}
	\newcommand{\vg}{\ensuremath{\mathrm{vol}_g}}
	\newcommand{\im}{\ensuremath{\mathrm{W}^{2,2}_{\iota}(\Sigma,N^n)}}
	\newcommand{\imm}{\ensuremath{\mathrm{W}^{2,2}_{\iota}(\Sigma,\R^3)}}
	\newcommand{\timm}[1]{\ensuremath{\mathrm{W}^{2,2}_{#1}(\Sigma,T\R^3)}}
	\newcommand{\tim}[1]{\ensuremath{\mathrm{W}^{2,2}_{#1}(\Sigma,TN^n)}}
	\renewcommand{\d}[1]{\ensuremath{\partial_{x_{#1}}}}
	\newcommand{\dg}{\ensuremath{\mathrm{div}_{g}}}
	\renewcommand{\Res}{\ensuremath{\mathrm{Res}}}
	\newcommand{\un}[2]{\ensuremath{\bigcup\limits_{#1}^{#2}}}
	\newcommand{\res}{\mathbin{\vrule height 1.6ex depth 0pt width
			0.13ex\vrule height 0.13ex depth 0pt width 1.3ex}}
	\newcommand{\ala}[5]{\ensuremath{e^{-6\lambda}\left(e^{2\lambda_{#1}}\alpha_{#2}^{#3}-\mu\alpha_{#2}^{#1}\right)\left\langle \nabla_{\vec{e}_{#4}}\vec{w},\vec{\mathbb{I}}_{#5}\right\rangle}}
	\setlength\boxtopsep{1pt}
	\setlength\boxbottomsep{1pt}
	\newcommand\norm[1]{%
		\setbox1\hbox{$#1$}%
		\setbox2\hbox{\addvbuffer{\usebox1}}%
		\stretchrel{\lvert}{\usebox2}\stretchrel*{\lvert}{\usebox2}%
	}
	\allowdisplaybreaks
	\newcommand*\mcup{\mathbin{\mathpalette\mcapinn\relax}}
	\newcommand*\mcapinn[2]{\vcenter{\hbox{$\mathsurround=0pt
				\ifx\displaystyle#1\textstyle\else#1\fi\bigcup$}}}
	\def\Xint#1{\mathchoice
		{\XXint\displaystyle\textstyle{#1}}%
		{\XXint\textstyle\scriptstyle{#1}}%
		{\XXint\scriptstyle\scriptscriptstyle{#1}}%
		{\XXint\scriptscriptstyle\scriptscriptstyle{#1}}%
		\!\int}
	\def\XXint#1#2#3{{\setbox0=\hbox{$#1{#2#3}{\int}$ }
			\vcenter{\hbox{$#2#3$ }}\kern-.58\wd0}}
	\def\ddashint{\Xint=}
	\newcommand{\dashint}[1]{\ensuremath{{\Xint-}_{\mkern-10mu #1}}}
	\newcommand\ccancel[1]{\renewcommand\CancelColor{\color{red}}\cancel{#1}}
	\newcommand\colorcancel[2]{\renewcommand\CancelColor{\color{#2}}\cancel{#1}}
		\renewcommand{\thetheorem}{\thesection.\arabic{theorem}}
		
	\section{Introduction}

	The classification of Willmore spheres in $S^3$ by Robert Bryant (\cite{bryant}) shows that any non-minimal Willmore sphere $\phi:\P^1\rightarrow S^3$ is the inverse stereographic projection of a complete minimal surface in $\R^3$ with embedded planar ends. Therefore, the classification of Willmore spheres in $S^3$ reduces to the classification of certain minimal surfaces in $\R^3$. It is relatively easy to find explicit examples of complete minimal surfaces with any even number $2d\geq 4$ of embedded planar ends (\cite{bryant}, \cite{kusnerpacific}, \cite{bryant3}), R. Bryant showed that there were such minimal surfaces with $3$, $5$ or $7$ ends.
	 
	Using the framework of the Klein correspondence (\cite{bryant3}, \cite{bryantcurves}), we show that there exists essentially a unique minimal surface with $8$ embedded flat ends, while there are no minimal surfaces with $9$ embedded planar ends. We also explain why the examples of the contrary which have been made in the literature are incorrect.
	
	
	
	\subsection{Motivation}
	
	We say that a non-constant meromorphic curve $f:\Sigma\rightarrow\C^3$ is a null curve if $\s{df}{df}=0$, or in a local coordinate $z:U\rightarrow\C$
	\begin{align*}
	\s{\p{z}f}{\p{z}f}=(\p{z}f_1)^2+(\p{z}f_2)^2+(\p{z}f_3)^2=0.
	\end{align*}
	The Weierstrass parametrisation shows that any non-planar minimal surface $\phi:\Sigma\rightarrow\R^3$ arises as $\phi=\Re(f)$, for such non-constant meromorphic null curve $f:\Sigma\rightarrow\C^3$, a condition which amounts as saying that $\phi$ is an immersion outside of the discrete set of poles of $f$. Furthermore, $\phi$ is a complete minimal surface with $n$ embedded planar ends if and only if $f$ has $n$ simple poles, \textit{i.e.} if there exists \emph{distinct} points $a_1,\cdots,a_n\in \C$ and vectors $v^0\in \C^3$, $v^1,\cdots,v^n\in \C^{3}\setminus\ens{0}$ such that
	\begin{align*}
	f(z)=v^0+\frac{v^1}{z-a_1}+\cdots+\frac{v^n}{z-a_n}.
	\end{align*}
	This is easy to see that $\s{f'(z)}{f'(z)}=0$ is an over-determined system in $(v^0,v^1,\cdots,v^n)$ and $(a_1,\cdots,a_n)$, and the direct approach for $n\geq 5$ seems quite difficult, so we will describe below another approach of R. Bryant using the \emph{Klein correspondence} on which the classification is based on.

	\subsection{The Klein correspondence}
		
		We adopt the notations of \cite{bryantcurves}, and we recall the Klein correspondence (see Proposition $4$ \cite{bryantcurves}). Let $\Sigma$ be a compact connected Riemann surface, and $f:\Sigma\rightarrow \P^n$ be a non-degenerate holomorphic curve,  \textit{i.e.} $f(\Sigma)$ is not contained in some hyperplane $H^{n-1}\subset \P^n$. Then this is known that $\mathrm{deg}(f)\geq n$ (see for example \cite{kollar}).		
		 It is possible to construct a well-defined (see \cite{griffiths}, 2.4) family of associated curves $f_k:\Sigma\rightarrow \P(\Lambda^k\C^{n+1})\simeq \mathbb{P}^{\binom{n+1}{k}-1}$ with $1\leq k\leq n$, such that for any local complex coordinate $z:U\rightarrow\C$ (where $U\subset \Sigma$ is an open set), 
		\begin{align*}
		f_k=[F\wedge \p{z}F\wedge\cdots\wedge \p{z}^{k-1}F]
		\end{align*}
		if $f=[F]$, where $F:U\rightarrow \C^{n+1}$ is a \emph{non-vanishing} holomorphic map. As we shall see, even if $f$ is non-degenerate, $f_k$ need not be non-degenerate.
		Notice that $f_k$ is the projectivization of Calabi's holomorphic form (see \cite{calabi})
		\begin{align*}
		F\wedge \partial F\wedge \cdots\wedge \partial^{k-1}F=F(z)\wedge \p{z}F(z)\wedge \cdots \wedge \p{z}^{k-1}F(z)\,dz^{\frac{k(k-1)}{2}},
		\end{align*}
		where in a local complex coordinate $z$, we have for all $1\leq j\leq k-1$
		\begin{align*}
		\partial^j F=\p{z}^jFdz^j.
		\end{align*}
		We now introduce the necessary definitions related to ramification divisors, and we adopt the same notations as \cite{bryantcurves} (see also \cite{griffiths} $1.4$).
		For all $p\in \Sigma$, there exists a basis $(v^0,\cdots,v^n)$ of $\C^{n+1}$ and holomorphic functions $h_0,\cdots,h_n$ on $\Sigma$ such that
		\begin{align*}
			f=\left[h_0\,v^0+h_1\,v^1+\cdots+h_n\,v^n\right]
		\end{align*} 
		and satisfying
		\begin{align}\label{hi}
			0=\mathrm{ord}_p(h_0)<\mathrm{ord}_p(h_1)<\cdots<\mathrm{ord}_p(h_n),
		\end{align}
		where $\mathrm{ord}_p(h_i)$ is the order of vanishing of $h_i$ at $p\in \Sigma$. We define the $i$-ramification degree of $f$ at $p$ for all $1\leq i\leq n$ by
		\begin{align*}
			r_i(f,p)=\mathrm{ord}_p(h_i)-\mathrm{ord}_p(h_{i-1})-1\geq 0,
		\end{align*}
		which is well-defined independently on the choice of basis $(v^0,\cdots,v^n)$ of $\C^{n+1}$ and of $\ens{h_i}_{1\leq i\leq n}$ satisfying \eqref{hi}.
		As for all but finitely many $p\in \Sigma$ do we have $r_i(f,p)>0$, we define for $1\leq i\leq n$ the $i$-th ramification divisor $R_i(f)$ on $\Sigma$ by 
		\begin{align*}
			R_i(f)=\sum_{p\in \Sigma}^{}r_i(f,p)\cdot p.
		\end{align*}
		We also mention the important set of relations between the ramification divisors of $f$ and of its associated curves:
		\begin{align*}
			R_i(f)=R_1(f_i),\quad R_i(f_n)=R_{n+1-i}(f),\quad 1\leq i\leq n.
		\end{align*}
		However, for $2\leq k\leq n$ and $2\leq i\leq \binom{n+1}{k}-1$, the branched divisor $R_i(f_k)$ cannot be computed solely with respect to $R_l(f)$ in general, with the notable exception of the contact curves as we shall see.
		
		\begin{defi}
			We say that $p\in \Sigma$ is a branch point of $f$ of order $\theta_0\geq 1$ if $r_1(f,p)=\theta_0$.  If $R_1(f)=0$, we say that $f$ is unbranched or equivalently an immersion in a neighbourhood of $p\in \Sigma$.
		\end{defi}
		
		Now assume that $n=3$. We say that a holomorphic curve $f:\Sigma\rightarrow\P^3$ is  a contact curve if there exists a symplectic form $\beta$ on $\C^4$ such that
		$f_2(S)\subset \P(\beta^{\perp})\subset \P(\Lambda^2\C^4)\simeq\P^5$. As up to linear transformation, we have (in the standard coordinates $(x_1,x_2,x_3,x_4)$ of $\C^4$)
		\begin{align*}
		\beta=dx_1\wedge dx_2+dx_3\wedge dx_4,
		\end{align*}
		we notice that $\beta^\perp$ is a $5$-dimensional sub-vector space of $\Lambda^2\C^4\simeq \C^6$.

		Now, we embed $\C^3\hookrightarrow \P^4$ as the null quadric $Q^3\subset \P^4$, also called \emph{Klein's quartic}, by the map 
		\begin{align*}
		x\mapsto [1,x,\s{x}{x}],
		\end{align*}	
		where the null quadric $Q^3\subset \P^4$ is defined by the homogeneous equations
		\begin{align}\label{quadric}
		X_0X_4-X_1^2-X_2^2-X_3^2=0.
		\end{align}
		Coming back to minimal surfaces, we mention the following correspondence with algebraic curves in $Q^3$.
		\begin{prop}[Bryant, \cite{bryantcurves}, Proposition $3$]\label{prop3}
			If $f:\Sigma\rightarrow \C^3$ is a meromorphic null curve with $d$ simple poles and no other poles, then the completed null curve $\widetilde{f}:\P^1\rightarrow Q^3\subset \P^4$ has degree $d$. Furthermore, if $f$ is an immersion outside of its poles, then $\widetilde{f}:\Sigma\rightarrow Q^3$ is also an immersion.
		\end{prop}
	
	    Coming back to the ramification divisors, we mention the remarkable fact that the ramifications divisors of the associate curve $f_2$ of a contact curve $f:\Sigma\rightarrow \P^3$ can be expressed solely with respect to the ramification divisors of $f$.
	    
	    \begin{prop}[Bryant, \cite{bryantcurves}, Proposition $1$]\label{nonlinear}
	    	Let $f:\Sigma\rightarrow \P^3$ a non-linear contact curve. Then $f$ is non-degenerate, and $f_2:\Sigma\rightarrow \P(\beta^{\perp})\simeq \P^4$ is non-degenerate as an algebraic curve in $\P^4$. Furthermore, we have
	    	\begin{align*}
	    		&R_1(f)=R_3(f)\\
	    		&R_1(f_2)=R_4(f_2)=R_2(f),\quad \text{and}\;\, R_2(f_2)=R_3(f_2)=R_1(f).
	    	\end{align*}
	    \end{prop}
		
		We finally come to the Klein correspondence (from \cite{klein}), which, in this strong form, is due to R. Bryant (\cite{bryantcurves}). 
		\begin{theorem}[Bryant, \cite{bryantcurves}, Propositions 1,2, and 4.]\label{kleinbryant}
			Let $f:\Sigma\rightarrow\P^3$ be a non-linear contact curve of associated symplectic form $\beta\in \Lambda^2((\C^4)^\ast)$, and $\Omega=\dfrac{1}{2}\beta^2$ be its associated volume form. Define the null hyperquadradric $Q^3$ of the scalar product $\s{\,\cdot\,}{\,\cdot\,}=\Omega(\,\cdot\,\wedge\,\cdot\,)$ by
			\begin{align*}
			Q^3=\P(\beta^\perp)\cap\ens{[v],\s{v}{v}=0}.
			\end{align*}
			Then $f_2:\Sigma\rightarrow \P(\Lambda^2\C^4)\simeq \P^5$ has image in $Q^3\subset \P(\beta^\perp)\simeq \P^4$, and is a non-degenerate null curve as a curve in $\P(\beta^\perp)\simeq \P^4$.
			
			Conversely, if $g:\Sigma\rightarrow Q^3\subset \P^4$ is a null curve whose image is not contained in a line $\P^1\subset \P^4$, then $g=f_2$ for a unique non-degenerate contact curve $f:\Sigma\rightarrow \P^3$.
			
		\end{theorem}

		Now, let us consider a contact holomorphic curve $f:\Sigma\rightarrow\P^3$ of symplectic form $\beta$, that is not contained in a line $\P^1\subset \P^3$, and let $f_2:\Sigma\rightarrow Q^3\subset \P(\beta^\perp)\simeq \P^4$ be its associated curve from Klein correspondence.
		The Plücker formulae and Proposition 1 of \cite{bryantcurves} (which proves that $R_1(f)=R_3(f)$, $R_1(f_2)=R_4(f_2)=R_2(f)$ and $R_2(f_2)=R_3(f_2)=R_1(f)$) show if $\Sigma$ has genus $g$ that
		\begin{align*}
		&4\deg(f)+12(g-1)=4r_1(f)+2r_2(f)\\
		&5\deg(f_2)+20(g-1)=5r_1(f)+5r_2(f),
		\end{align*}
		where $r_i(f)=\deg R_i(f)$.
		Notice that this implies that $r_2(f)$ is even.
		In particular, if $\Sigma$ has genus $0$, then
		\begin{align}\label{klein}
		&\deg(f)=3+r_1(f)+\frac{1}{2}r_2(f)\nonumber\\
		&\deg(f_2)=4+r_1(f)+r_2(f).
		\end{align}
		
		\begin{defi}
			We say that a non-degenerate contact curve $f:\Sigma\rightarrow \P^3$ is totally ramified if $r_1(f)$ is maximal, that is $r_2(f)=0$ and $r_1(f)=\mathrm{deg}(f)+3(g-1)$.
		\end{defi}
		
		Now, recall that one of the main results of \cite{bryant3} or \cite{bryantcurves} is to show the following theorem, which is equivalent to the non-existence of complete minimal surfaces with $5$ or $7$ embedded planar ends in $\R^3$. We can also easily check that for $2$ or $3$ ends, there is no such objects by a direct algebraic computation from the Weierstrass parametrisation. For example, the only complete minimal surface with $2$ embedded ends is the catenoid (see \cite{schoenPlanar}), whose ends are not planar. Furthermore, $3$ ends are excluded as the corresponding contact curve $f:\P^1\rightarrow \P^3$ would have degree $2$   
		
		\begin{theorem}[\label{deg7}Bryant, \cite{bryant3}, \cite{bryantcurves}]\label{main0}
			Unbranched non-linear null curves $g:\P^1\rightarrow Q^3\subset \P^4$ cannot have degree $5$ or $7$.
		\end{theorem}
		
		Using the link with minimal surfaces of Proposition \ref{prop3} and Theorem \ref{deg7}, we obtain the following non-existence result concerning minimal surfaces in $3$-space.
		
		\begin{theorem}[Bryant \cite{bryant3}, \cite{bryantcurves}]\label{mainb}
			Complete minimal surfaces of genus $0$ in $\R^3$ with  $3$, $5$, or $7$ embedded flat ends (and no other ends) do not exist.
	\end{theorem}

\subsection{Unbranched null curves of even degree}

The cases $d=4$, $d=6$ were completely classified by R. Bryant (see \cite{bryant}, \cite{bryant3}), and examples for even $d\geq 8$ were given by R. Kusner (see \cite{kusnerpacific}). Furthermore, there is a simple example given in \cite{bryant3} of curves of degree $2d\geq 4$ as the associate curve of the non-degenerate contact curve $f_d:\P^1\rightarrow\P^4$ defined by 
\begin{align}\label{fd}
f_d=\left[-\left(\frac{1}{2d-1}\right)v^0+z^{d-1}v^1+z^dv^2+z^{2d-1}v^3\right]
`\end{align}
It is a contact curve for the non-degenerate symplectic structure
\begin{align}\label{rel2}
\beta=\xi_0\wedge \xi_3+\xi_1\wedge \xi_2,
\end{align}
where $(\xi_0,\xi_1,\xi_2,\xi_3)$ is the dual base of the base $(v^0,v^1,v^2,v^3)$ of $\C^4$. Its branching divisor is
\begin{align*}
R_1(f)=(d-2)p+(d-2)q
\end{align*}
if $p,q\in\P^1$ correspond to the zero and the pole of the standard meromorphic coordinate $z$ on $\P^1=\C\cup\ens{\infty}$.
Indeed, we easily compute that
\begin{align*}
(f_d)_2=\left[-\frac{d-1}{2d-1}\,v^0\wedge v^1-\frac{d}{2d-1}z\,v^0\wedge v^2+z^{d}(v^1\wedge v^2-v^0\wedge v^3)+dz^{2d-1}\,v^1\wedge v^3+(d-1)z^{2d}\,v^2\wedge v^3\right]
\end{align*}
and as $2d>2d-1>d>1$ for all $d\geq 2$, we deduce that $g$ is linearly full in the projectivization of $W=\beta^\perp$, where $\beta$ is the non-degenerate symplectic form on $\C^4$ given by \eqref{rel2}. We will see that $f_d$ is up to projective equivalence the only totally ramified non-degenerate contact curve of degree $2d-1$ (at least for $d=4$).

\subsection{Statement of the results}

\setcounter{theorem}{0}

\renewcommand{\thetheorem}{\Alph{theorem}}

The first result permits to classify contact curves whose dual is an immersion in Klein's quadric.
	
	\begin{theorem}\label{TA}
		Let $f:\P^1\rightarrow \P^3$ be a non-degenerate totally ramified contact curve of degree $d\leq 9$. Then $d$ is odd. 
	\end{theorem}	

Using the Klein correspondence, this result permits to generalise Proposition $3$ of \cite{bryant3}. 

\begin{theorem}\label{TB}
	If $g:\P^1\rightarrow Q^3\subset \P^4$ is an unbranched holomorphic null immersion of degree $4\leq d\leq 9$, then $d$ is \emph{even} and $g$ is equivalent to the dual curve of $f_{d-1}:\P^1\rightarrow \P^3$ up to re-parametrisation in $\P^1$ and the action of the holomorphic automorphism $\mathrm{SO}(5,\C)$ of $Q^3$.  
\end{theorem}

By the Klein correspondence (\cite{bryant3}), we deduce the following result. 


 \begin{cor}
 	Complete minimal surfaces of genus $0$ in $\R^3$ with exactly $9$ embedded planar ends (and no other ends) do not exist. 
 \end{cor}
 It seems likely that this result holds for every odd number at least $11$ (we show that in general many branched divisors are excluded, see Sections \ref{excluded1}, \ref{excluded2}).
 \begin{conjecture}
 	Let $\phi:S^2\setminus\ens{p_1,\cdots,p_d}\rightarrow \R^3$ be a non-planar minimal surface with embedded planar ends. Then $d\geq 4$ is even.  
 \end{conjecture}
 
 

\begin{remast}
	If this conjecture held, it would draw a remarkable parallel between Willmore surfaces $S^2\rightarrow \R^3$ and harmonic maps $\R^3\rightarrow S^2$. Indeed, it would imply that for all non-completely umbilic Willmore immersion $\phi:S^2\rightarrow \R^3$,
	\begin{align*}
		W(\phi)\in 8\pi \N,
	\end{align*}
	while for all variational Willmore sphere $\phi:S^2\rightarrow \R^3$ (thanks to the combined results of \cite{classification}, \cite{sagepaper}, \cite{blow-up}, \cite{blow-up2}), we have
	\begin{align*}
		W(\phi)\in 4\pi \N.
	\end{align*}
    This is reminiscent of the work of Lin and Rivi\`{e}re (\cite{linriv}) on the energy quantization of harmonic maps, where they show in particular that for stationary harmonic maps $u:\R^3\rightarrow S^2$, we have
	\begin{align*}
		\lim\limits_{r\rightarrow \infty}\frac{1}{r}\int_{B(0,r)}\frac{1}{2}|\D u|^2dx=4\pi\, d
  	\end{align*}
	for some $d\in \N$ and that if $u:\R^3\rightarrow S^2$ is \emph{smooth} we have the stronger
	\begin{align*}
		\lim\limits_{r\rightarrow \infty}\frac{1}{r}\int_{B(0,r)}\frac{1}{2}|\D u|^2dx=8\pi\, d
	\end{align*}
	for some $d\in \N$.
\end{remast}



	
\section{Higher genus minimal surfaces with flat ends}

Finally, one can wonder what happens in higher genus. To our knowledge, the only known examples are in genus $1$. First, recall that the Jorge-Meeks formula (\cite{jorge}) shows that for any complete minimal immersion $\phi:\Sigma\setminus\ens{p_1,\cdots,p_d}\rightarrow \R^3$ of finite total curvature, if $p_1,\cdots,p_d$ have respective multiplicities $m_1,\cdots,m_d\geq 1$, then
\begin{align*}
\int_{\Sigma}K_g\,d\vg=-4\pi\left(\gamma-1+\frac{1}{2}\sum_{j=1}^{d}\left(m_j+1\right)\right),
\end{align*} 
where $\gamma\in \N$ is the genus of $\Sigma$.
In particular, if $\Sigma$ has genus $1$, and $\phi$ has $d$ embedded ends, we find
\begin{align*}
\int_{\Sigma}K_g\,d\vg=-4\pi d.
\end{align*}
The case $d=1$ is impossible, as the only complete minimal surfaces with total curvature $-4\pi$ are the catenoid and the Enneper surface (\cite{osserman}). The case $d=2$ is also impossible by the uniqueness of the catenoid as the only complete minimal surface with two embedded ends. The case $d=3$ is impossible by an argument of R. Kusner and N. Schmitt. Finally, C. Costa gave an example with $4$ ends (\cite{costa2}), and Kusner-Schmitt computed the moduli space of these minimal tori with $4$ flat ends. Additional examples of any \emph{oven} number of ends (at least $6$) were provided by E. Shamaev (\cite{shamaev}). We find it of interest that the only known examples have an even number of ends, and a bold conjecture might be to say that examples with an odd number of flat ends do not exist. Furthermore, it seems plausible that complete minimal surfaces of arbitrary genus with an even number of flat ends exist, by "adding handles" to the minimal surfaces constructed by R. Bryant and R. Kusner (see for example \cite{helicoid} for the construction of prescribed genus helicoids). 



\renewcommand{\thetheorem}{\thesection.\arabic{theorem}}

\setcounter{theorem}{0}

\section{Impossible divisors for unbranched null immersions}\label{excluded1}

We have already seen that for $d=1,2,3,4$, there is no unbranched null curve $g:\P^1\rightarrow Q_3\subset \P^4$ of degree $2d+1$. Let $d\geq 2$, and we suppose that there exists an unbranched null curve $g:\P^1\rightarrow Q_3\subset \P^4$, and we let $f:\P^1\rightarrow \P^3$ the associate contact curve from the Klein correspondence. As $g$ is unbranched, $f$ has degree $2d\geq 4$, and we have by the Pl\"{u}cker formulae
\begin{align*}
r_1(f)=2d-3\geq 1.
\end{align*}
Now, suppose that for some $p\in \P^1$, we have $R_1(f)=(2d-3)\cdot p$. Taking $p=0$, we see that for some $\lambda_1,\cdots, \lambda_{2d-3}\in \C$, we have and for some vectors $(v^0,v^1,v^2,v^3)\in \C^4$
\begin{align*}
f=\left[(1+\lambda_1z+\lambda_2z^2+\cdots +\lambda_{2d-3}z^{2d-3})v^0+z^{2d-2}v^1+z^{2d-1}v^2+z^{2d}v^3\right].
\end{align*}
As $f$ is non-degenerate, $(v^0,v^1,v^2,v^3)$ must be a base of $\C^4$. Now, we compute
\begin{align*}
&f_2=\Big[(2d-2)v^0\wedge v^1z^{2d-3}+((2d-1)v^0\wedge v^2+(\ast\ast\ast)v^0\wedge v^1)z^{2d-2}\\
&+(2dv^0\wedge v^3+(\ast\ast\ast)v^0\wedge v^1+(\ast\ast\ast)v^0\wedge v^2)z^{2d-1}+(\ast\ast\ast)z^{2d}+\cdots +(\ast\ast\ast)z^{4d-5}\\
&+(3\lambda_{2d-3}v^0\wedge v^3+v^1\wedge v^2)z^{4d-4}
+(2v^1\wedge v^3)z^{4d-3}+(v^2\wedge v^3)z^{4d-2}\Big].
\end{align*}
Looking at first three and the last three lines, we see that $f_2$ is linearly full in $\P(\Lambda^2(\C^4))$, so $f$ cannot be a contact curve of degree $2d\geq 4$, as $4d-4>2d-1$.

\begin{rem}
	Notice that the same proof would work for degree $2d\geq 4$ unbranched null curves in Klein's quadric $Q^3\subset \P^4$.
\end{rem}

Actually, we can also obtain this directly thanks of the following lemma.

\begin{lemme}\label{cont}
	Let $f:\P^1\rightarrow \P^3$ be a completely ramified contact curve of degree $d\geq 5$. Then for all $p\in \P^1$, we have
	\begin{align*}
		r_1(f,p)\leq  \frac{d-3}{2}.
	\end{align*}
	and in particular
	\begin{align*}
		\P^1\cap\ens{p:r_1(f,p)>0}\geq 2.
	\end{align*}
\end{lemme}
\begin{proof}
	Assume that $p\in \P^1$ is such that $2r_1(f)>d-3$. Then if $z$ is the standard meromorphic coordinate of $\P^1$, we can assume that $p$ corresponds to the zero of $z$ and there exists $(v^0,v^1,\cdots,v^d)\in \C^4$ such that
	\begin{align*}
		f=\left[v^0+zv^1+\cdots z^dv^d\right].
	\end{align*}
	Now, if $a=r_1(f,p)\geq 1$, $v^1,\cdots,v^a$ must be multiples of $v^0$ and we obtain for some $\lambda_1,\cdots,\lambda_a\in \C$
	\begin{align*}
		f=\left[\left(1+\lambda_1z+\cdots+\lambda_az^a\right)v^0+z^{a+1}v^{a+1}+z^{a+2}v^{a+2}+\cdots+z^dv^d\right].
	\end{align*}
	As $r_2(f,p)=0$, we see that $(v^0,v^{a+1},v^{a+2})$ is free and as $r_3(f,p)=r_1(f,p)$, we must have for some basis $(w^0,w^1,w^2,w^3)$ of $\C^4$ the expansion
	\begin{align*}
	f=\left[\left(1+\cdots\right)v^0+\left(z^{a+1}+\cdots\right)w^1+\left(z^{a+2}+\cdots\right)w^2+\left(z^{2a+3}+\cdots\right)w^3\right].
	\end{align*}
	so we have $v^j\in \mathrm{Span}(v^0,v^{a+1},v^{a+2})$ for all $a+3\leq j\leq  2a+2$.
	However, as $d\leq 2a+2$, we have $v^j\in \mathrm{Span}(v^0,v^{a+1},v^{a+2})$ for all $a+3\leq j\leq d$ so $f(\P^1)\subset \P(\mathrm{Span}(v^0,v^{a+1},v^{a+2}))\simeq \P^2\subset \P^3$ so $f$ is degenerate, contradiction.
\end{proof}

Therefore, we obtain the following partial result.

\begin{prop}
	Suppose that there exists an unbranched null curve $g:\P^1\rightarrow Q_3\subset\P^4$ of degree $2d+1$, and let $f:\P^1\rightarrow\P^3$ be the associated contact curve given by the Klein correspondence. Then the ramification divisor of $f$ consists of at least three distinct points with multiplicity.
\end{prop}
\begin{proof}
	We have $\mathrm{deg}(f)=2d$ and $r_1(f)=2d-3$, and by Lemma \ref{cont}, we have for all $p\in \P^1$
	\begin{align*}
		r_1(f,p)\leq \frac{2d-3}{2}=d-\frac{3}{2}
	\end{align*}
	which implies that $r_1(f,p)\leq d-2$ and as $2(d-2)<2d-3=r_1(f)$, the algebraic curve $f$ must have at least three distinct branch points.
\end{proof}

We can refine this result thanks of the following lemmas, using Proposition $1$ of \cite{bryantcurves}.

\begin{lemme}\label{abel}
	Let $f:\P^1\rightarrow \P^3$  be a non-degenerate totally ramified contact curve of degree $d\geq 4$ (\textit{i.e.} $r_1(f)=d-3$). Then for all $p\in \P^1$ such that $2r_1(f,p)+3<d$, and for all $q\in \P^1\setminus\ens{p}$, we have the estimate
	\begin{align}\label{abel0}
	2\,r_1(f,p)+r_2(f,q)\leq d-4.
	\end{align}
	In particular, if $f$ has even degree, then \eqref{abel0} holds for all $p\neq q$.
\end{lemme}
\begin{proof}
	Let $p,q\in \P^1$ two distinct points and let $a=r_1(p,f)$, $b=r_2(f,q)$. Then we can assume that $p$ corresponds to the zero of the meromorphic coordinate $z$ and that $q$ corresponds to its pole. Then for some vectors $v^0,\cdots,v^d\in \C^4$ spanning $\C^4$, and some scalars $\lambda_j,\nu_j\in \C$, we can write $f$ as
	\begin{align}\label{comp1}
	f=\left[\left(1+\lambda_1z+\cdots+\lambda_az^a\right)v^0+z^{a+1}v^{a+1}+z^{a+2}v^{a+2}+\cdots+\left(\pi_{d-b}z^{d-b}+\cdots+z^d\right)v^d\right].
	\end{align}
	Now, assuming that $2d+3<d$, we can make a linear change of variable such that $\pi_{2a+3}=0$ (up to modifying the other coefficients).
	As $R_2(f)=0$ and $R_3(f)=R_1(f)$, we can write $f$ as
	\begin{align}\label{comp2}
	f=\left[\left(1+\cdots\right)w^0+\left(z^a+\cdots\right)w^1+\left(z^{a+1}+\cdots\right)w^2+\left(z^{2a+3}+\cdots\right)w^3\right]
	\end{align}
	for some basis $(w^0,w^1,w^2,w^3)$ of $\C^4$. In particular, assuming that $2d+3\geq d-b$, this implies that $v^{j}\in \mathrm{Span}(v^0,v^{a+1},v^{a+2})$ for all $a+3\leq j\leq d-b-1$, so for some scalars $\lambda_j,\mu_j,\nu_j\in \C$, we have
	\begin{align*}
	f&=\Big[\left(1+\lambda_1z+\cdots+\lambda_az^a+\lambda_{a+3}z^{a+3}+\cdots+\lambda_{d-b-1}z^{d-b-1}\right)v^0+\big(z^{a+1}+\mu_{a+3}z^{a+3}+\cdots\\ &+\mu_{d-b-1}z^{d-b-1}\big)v^{a+1}
	+\left(z^{a+2}+\nu_{a+3}z^{a+3}+\cdots+\nu_{d-b-1}z^{d-b-1}\right)v^{a+2}+\left(\pi_{d-b}z^{d-b}+\cdots+z^d\right)v^d \Big].
	\end{align*}
	In particular we see that $(v^0,v^{a+1},v^{a+2},v^{d})$ is a basis of $\C^4$. 
	Furthermore, by \eqref{comp2}, we must have $\pi_j=0$ for all $d-b\leq j\leq 2a+2$ (otherwise we would have $r_3(f,p)<a$), so recalling that we also have $\pi_{2a+3}=0$ we finally obtain
	\begin{align*}
	f&=\Big[\left(1+\lambda_1z+\cdots+\lambda_{d-b-1}z^{d-b-1}\right)v^0+\big(z^{a+1}+\mu_{a+3}z^{a+3}+\cdots
	+ \mu_{d-b-1}z^{d-b-1}\big)v^{a+1}\\
	&+\left(z^{a+2}+\nu_{a+3}z^{a+3}+\cdots+\nu_{d-b-1}z^{d-b-1}\right)v^{a+2}+\left(\pi_{2a+4}z^{2a+4}+\cdots+z^d\right)v^d \Big],
	\end{align*}
	which contradicts \eqref{comp2}, as this expression shows that $r_3(f,p)>a=r_1(f,p)=r_3(f,p)$, as the ramifications divisors $R_1,R_2,R_3$ are independent of coordinates and of the choice of the base of $\C^4$ and of the meromorphic function realising the appropriate Taylor expansion.
\end{proof}
\begin{cor}\label{cor}
	Let $f:\P^1\rightarrow \P^3$ be a non-degenerate and totally ramified contact curve of degree $2d\geq 4$. Then for all $p\in \P^1$, there holds
	\begin{align*}
	r_1(f,p)\leq d-3.
	\end{align*}
\end{cor}
\begin{proof}
	As $R_1(f)\neq (d-3)\cdot p$ for some $p\in \P^1$, there exists $q\in \P^1\setminus\ens{p}$ such that $r_1(f,q)\geq 1$, so by \eqref{abel0}
	\begin{align*}
	2\,r_1(f,p)+1\leq 2r_1(f,p)+r_2(f,p)\leq 2d-4
	\end{align*}
	so that $2\,r_1(p,f)\leq 2d-5$, implying the claim.
\end{proof}
\begin{prop}
	Let $f:\P^1\rightarrow \P^3$ be a non-degenerate and totally ramified contact curve of degree $2d$. Then we have
	\begin{align*}
	\mathrm{card}\;\left(\P^1\cap\ens{p:r_1(f,p)>0}\right)\geq 4.
	\end{align*}
\end{prop}
\begin{proof}
	As by the previous Corollary \ref{cor}, we have $r_1(f,p)\leq d-3$ for all $p\in \P^1$, we deduce in particular that for all $p,q\in \P^1$
	\begin{align*}
	r_1(f,p)+r_1(f,q)\leq 2d-6<2d-3=r_1(f)
	\end{align*}
	so there exists at least three distinct points $p,q,r\in \P^1$ such that $r_1(f,p)>0$ (\textit{i.e.} $f$ has at least three distinct branch points).
	Notice that the bound $r_1(f,p)\leq d-2$ would suffice for this argument.
	Now, suppose that $f$ has exactly three distinct branch points $p,q,r\in \P^1$ of respective multiplicities $a\geq b\geq c\geq 1$. As
	\begin{align*}
	r_1(f)=a+b+c=2d-3,
	\end{align*}
	we deduce that $a\geq \dfrac{2d}{3}-1$. Now, by Lemma \ref{abel}, we have
	\begin{align*}
	c\leq b\leq 2d-4-2a\leq \frac{2d}{3}-2,
	\end{align*}
	so that
	\begin{align*}
	2d-3=a+b+c\leq a+2(2d-4-2a)=-3a+4d-8,
	\end{align*}
	so that $3a\leq 2d-5$, and   this implies as $a\geq b\geq c$ that
	\begin{align*}
	2d-3=a+b+c\leq 3a\leq 2d-5<2d-3,
	\end{align*}
	a contradiction.
\end{proof}

In particular, we recover without computations Theorem \ref{mainb}.

\begin{prop}\label{lazy}
	Let $f:\P^1\rightarrow\P^3$ a totally ramified contact curve of degree $2d\geq 4$. Then $2d\geq 8$, and if $d=8$, then $R_1(f)=p+q+r+s+t$ for five distinct points $p,q,r,s,t\in \P^1$.
\end{prop}
\begin{proof}
	If $f$ has degree $2d\geq 4$, then $r_1(f)=2d-3\geq 1$, and as $r_1(f,p)\leq d-3$ for all $p\in \P^1$, if $d\leq 3$, then $r_1(f,p)=0$ for all $p\in \P^1$, while $r_1(f)>0$, a contradiction. If $\mathrm{deg}(f)=2d=8$, so $r_1(f)=2d-3=5$, and $r_1(f,p)\leq d-3=1$ for all $p\in \P^1$ implies the claim.
\end{proof}

\section{Contact curves of odd degree}

It seems that there always exists a unique (up to the re-parametrisation in $\P^1$ and the action of the holomorphic automorphism $\mathrm{SO}(5,\C)$ of $Q^3$) unbranched null curves of even degree. For now, we only have the following very partial result.

\begin{prop}\label{p5}
	Let $f:\P^1\rightarrow \P^3$ is a non-degenerate contact curve of degree $2d-1\geq 3$ such that
	\begin{align*}
	\mathrm{card}\left(\P^1\cap\ens{p:r_1(f,p)>0}\right)\leq 2.
	\end{align*}
	Then $f=f_{d}$ given by \eqref{fd}.
\end{prop}
\begin{proof}
	For $d=2$, $\mathrm{deg}(f)=3$, so $f$ is a rational normal curve, so there is nothing to do. If $d\geq 3$, by Lemma \ref{cont}, $r_1(p,f)\leq d-2$ and $0<d-2<2d-4$, so $R_1(f)=(2d-4)\cdot p$ for some $p\in \P^1$. Therefore, $R_1(f)=m\cdot p+n\cdot q$ for two distinct points $p,q\in \P^1$. As $m,n\leq d-2$ and $m+n=2d-4$, we have $m=n=d-2$, so there exists a basis $(v^0,v^1,v^2,v^3)$ and $\lambda_j,\pi_j\in \C$
	\begin{align*}
	f=\left[\left(1+\lambda_1z+\cdots \lambda_{d-2}z^{d-2}\right)v^0+z^{d-1}v^1+z^{d}v^2+\left(\pi_1z^{d+1}+\cdots+\pi_{d-2}z^{2d-2}+z^{2d-1}\right)v^3\right]
	\end{align*}
	As $R_3(f)=R_1(f)$, we see that we must have $\lambda_j=\pi_k=0$ for all $1\leq j\leq d-2$ and $1\leq k\leq d-2   $, so $f=f_{d}$.
\end{proof}

\begin{theorem}\label{deg8}
	If $g:\P^1\rightarrow Q^3\subset \P^4$ is an unbranched holomorphic null immersion of degree $8$, then $g$ is equivalent to the dual curve of $f_{4}:\P^1\rightarrow \P^3$ up to re-parametrisation in $\P^1$ and the action of the holomorphic automorphism $\mathrm{SO}(5,\C)$ of $Q^3$.
\end{theorem}
\begin{proof}
	Let $f:\P^1\rightarrow \P^3$ a contact curve of degree $d$, then recall that
	\begin{align}\label{d-3}
	    &r_1(f)=d-3\nonumber\\
		&r_1(f,p)\leq \frac{d-3}{2}\;\,\text{for all}\;\, p\in \P^1.
	\end{align}
	Therefore, if $f:\P^1\rightarrow \P^3$ has degree $7$, $r_1(f)=4$, $r_1(f,p)\leq 2$ for all $p\in \P^1$, so the possible divisors are
	\begin{align*}
		&2\cdot p+2\cdot q\\
		&2\cdot p+q+r\\
		&p+q+r+s
	\end{align*}
	for some distinct $p,q,r,s\in  \P^1$. As we have already seen, the first one corresponds to $f_4$. Now, if $R_1(f)=2\cdot p+q+r$ or $R_1(f)=p+q+r+s$, then we have at least two branch points of order $1$, so we can assume that they corresponds to $z=0$ and $z=\infty$, so that
	\begin{align*}
		f=\left[(1+\lambda_1z)v^0+z^2v^1+z^3v^2+z^4v^3+z^5v^4+(\pi_1z^6+z^7)v^5\right]=[F(z)]
	\end{align*}
	Using $R_3(f)=R_1(f)$ and $R_2(f)=0$, we see that $(v^0,v^1,v^2,v^4)$ and $(v^5,v^4,v^3,v^1)$ are a basis of $\C^4$. Now, we can make a change of basis so that $(v^0,v^1,v^2,v^4)$ is orthogonal, by introducing some $\lambda_j,\mu_j,\nu_j\in \C$ such that
	\begin{align}\label{gram}
	   f&=\Big[\left(1+\lambda_1z+\lambda_2z^2+\lambda_3z^3+\lambda_5z^5\right)v^0+(z^2+\mu_1z^3+\mu_2z^5)v^1+(z^3+\nu_2z^5)v^2+z^4v^3+z^5v^4\nonumber\\
	   &+(\pi_1z^6+z^7)v^5\Big].
	\end{align}
	Indeed, we define recursively by the Gram-Schmidt orthogonalisation process (and denoting by $\s{\,\cdot\,}{\,\cdot\,}$ the Hermitian form of $\C^4$)
	\begin{align*}
		&\zeta_0=\left|v^0\right|, \quad \widetilde{v}^0=\frac{v^0}{|v^0|}\\
		&\zeta_1=\left|v^1-\s{v^1}{\widetilde{v}^0}\widetilde{v}^0\right|,\quad \widetilde{v}^1=\frac{v^1-\s{v^1}{\widetilde{v}^0}\widetilde{v}^0}{|v^1-\s{v^1}{\widetilde{v}^0}v^0|}\\
		&\zeta_2=\left|v^2-\s{v^2}{\widetilde{v}^0}\widetilde{v}^0-\s{v^2}{\widetilde{v}^1}\widetilde{v}^1\right|,\quad \widetilde{v}^2=\frac{v^2-\s{v^2}{\widetilde{v}^0}\widetilde{v}^0-\s{v^2}{\widetilde{v}^1}\widetilde{v}^1}{|v^2-\s{v^2}{\widetilde{v}^0}\widetilde{v}^0-\s{v^2}{\widetilde{v}^1}\widetilde{v}^1|}\\
		&\zeta_4=\left|v^4-\s{v^4}{\widetilde{v}^0}\widetilde{v}^0-\s{v^4}{\widetilde{v}^1}\widetilde{v}^1-\s{v^4}{\widetilde{v}^2}\widetilde{v}^2\right|,\quad \widetilde{v}^4=\frac{v^4-\s{v^4}{\widetilde{v}^0}\widetilde{v}^0-\s{v^4}{\widetilde{v}^1}\widetilde{v}^1-\s{v^4}{\widetilde{v}^2}\widetilde{v}^2}{|v^4-\s{v^4}{\widetilde{v}^0}\widetilde{v}^0-\s{v^4}{\widetilde{v}^1}\widetilde{v}^1-\s{v^4}{\widetilde{v}^2}\widetilde{v}^2|}.
	\end{align*}
	Therefore, we have
	\begin{align*}
		F(z)&=(1+\lambda_1z)v^0+z^2v^1+z^3v^2+z^4v^3+z^5v^4+(\pi_1z^6+z^7)v^5\\
		&=\zeta_0\left(1+\lambda_1\right)\widetilde{v}^0+z^2\left(\zeta_1\widetilde{v}^1+\s{v^1}{\widetilde{v}^0}\widetilde{v}^0\right)+z^3\left(\zeta_2\widetilde{v}^2+\s{v^2}{\widetilde{v}^0}\widetilde{v}^0+\s{v^2}{\widetilde{v}^1}\widetilde{v}^1\right)+z^4v^3\\
		&+z^5\left(\zeta_4\widetilde{v}^4+\s{v^4}{\widetilde{v}^0}\widetilde{v}^0+\s{v^4}{\widetilde{v}^1}\widetilde{v}^1+\s{v^4}{\widetilde{v}^2}\widetilde{v}^2\right)+\left(\pi_1z^6+z^7\right)v^5\\
		&=\left(\zeta_0+\zeta_0\lambda_1+\s{v^1}{\widetilde{v}^0}z^2+\s{v^2}{\widetilde{v}^0}z^3+\s{v^4}{\widetilde{v}^0}z^5\right)\widetilde{v}^0+\left(\zeta_1+\s{v^2}{\widetilde{v}^1}z^3+\s{v^4}{\widetilde{v}^1}z^5\right)\widetilde{v}^1\\
		&+\left(\zeta_2z^3+\s{v^4}{\widetilde{v}^2}z^5\right)\widetilde{v}^2+z^4v^3+\zeta_4z^5\widetilde{v}^4+\left(\pi_1z^6+z^7\right)v^5,
	\end{align*}
	and $(\widetilde{v}^0,\widetilde{v}^1,\widetilde{v}^2,\widetilde{v}^4)$ is an orthonormal basis of $\C^4$, so up to renaming and scaling each of these coefficients by a non-zero real constant, we can suppose that $f$ is given by \eqref{gram}, where $(v^0,v^1,v^2,v^4)$ is an \emph{orthogonal} basis of $\C^4$ (though not orthonormal in general).
	
	Therefore, as $(v^0,v^1,v^2,v^4)$ is an orthogonal basis of $\C^4$, and recalling that  $(v^1,v^3,v^4,v^5)$ is also a basis of $\C^4$ we must have
	\begin{align*}
		\mathrm{Span}(v^0,v^2)=\mathrm{Span}(v^3,v^5)\simeq \C^4,
	\end{align*}
	so there exists $\lambda_4,\lambda_6,\nu_1,\mu_3\in \C$ such that
	\begin{align*}
		\left\{\begin{alignedat}{1}
		&v^3=\lambda_4v^0+\nu_1v^2\\
		&v^5=\lambda_6v^0+\nu_3v^2
		\end{alignedat}\right.
	\end{align*}
	Therefore, we have up to renaming $v^4$ in $v^3$
	\begin{align*}
		F(z)&=\left(1+\lambda_1z+\lambda_2z^2+\lambda_3z^3+\lambda_4z^4+\lambda_5z^5+\lambda_6\left(\pi_1z^6+z^7\right)\right)v^0\\
		&+\left(z^2+\mu_1z^3+\mu_2z^5\right)v^1+\left(z^3+\nu_1z^4+\nu_2z^5+\nu_3\left(\pi_1z^6+z^7\right)\right)v^2+z^5v^3.
	\end{align*}
	so that
	\begin{align*}
		f_2=&\bigg[(2+\cdots)v^0\wedge v^1+(3z+\cdots)v^0\wedge v^2+(5z^3+\cdots)v^0\wedge v^3+(z^3+\cdots)v^1\wedge v^2+z^5\left(3+2\mu_1z\right)v^1\wedge v^3\\
		&+(2z^6+\cdots)v^2\wedge v^3\bigg]=[F_2(z)].
	\end{align*}
	However, the coefficient $z^5\left(3+2\mu_1z\right)$ in $v^1\wedge v^3$ has at most one zero on $\C\setminus\ens{0}$ and as $f$ has two branch points (with multiplicity) outside on $\C\setminus\ens{0}$, $F_2$ must have two zeroes (with multiplicity) on $\C\setminus\ens{0}$, which is a contradiction.
\end{proof}

\section{Unbranched null immersion of degree $9$}

\begin{theorem}\label{deg9}
	An unbranched null curve $g:\P^1\rightarrow Q^3\subset \P^4$ cannot be of degree $9$.
\end{theorem}
\begin{proof}
Thanks of Proposition \ref{lazy}, the associate contact curve $f:\P^1\rightarrow \P^3$ of an hypothetical degree $9$ unbranched null curve $g:\P^1\rightarrow Q^3\subset \P^4$ has degree $8$ and is such that $R_1(f)=p+q+r+s+t$ for some distinct points $p,q,r,s,t\in \P^1$.

 Then we can assume that $0$ and $\infty$ are branched points of order $1$, so that
\begin{align}\label{eq2}
	f=\left[(1+\lambda_1z)v^0+z^2v^1+z^3v^2+z^4v^3+z^5v^4+z^6v^5+(\pi_1z^7+z^8)v^6\right]
\end{align}
As $R_2(f)=0$, and $R_3(f)=R_1(f)$, and as $f$ has a branch point of order $1$ at $0$
there exists a basis $(w^0,w^1,w^2,w^3)$ of $\C^4$, such that
\begin{align*}
	f=\left[(1+\cdots)w^0+(z^2+\cdots)w^1+(z^3+\cdots)w^2+(z^5+\cdots)w^3\right]
\end{align*}
where $+\cdots$ designs terms of higher order, looking at \eqref{eq2}, we see that $(v^0,v^1,v^2,v^4)$ is also a basis of $\C^4$.

Likewise, as $f$ has a branch point of order $1$ at $z=\infty$, the family $(v^6,v^5,v^4,v^2)$ must be a basis of $\C^4$.
 As the two families $(v^0,v^1,v^2,v^4)$ and $(v^2,v^4,v^5,v^6)$ are a basis of $\C^4$. Now, by a change of basis, we can assume than $(v^0,v^1,v^2,v^4)$ is an orthonormal basis, if
 \begin{align*}
 	F(z)&=(1+\lambda_1z+\lambda_2z^2+\lambda_3z^3+\lambda_5z^5)v^0+(z^2+\mu_1z^3+\mu_3z^5)v^1+(z^3+\nu_2z^5)v^2+z^4v^3+z^5v^4+z^6v^5\\
 	&+(\pi_1z^7+z^8)v^6.
 \end{align*}
 As $(v^2,v^4,v^5,v^6)$ is a basis of $\C^4$, we have $v^0,v^1\in \mathrm{Span}(v^2,v^4,v^5,v^6)$ but as $(v^0,v^1,v^2,v^4)$ is orthonormal, we must actually have $v^0,v^1\in \mathrm{Span}(v^5,v^6)$ and as $\mathrm{Span}(v^0,v^1)\simeq \C^2$, we finally deduce that 
\begin{align*}
	\mathrm{Span}(v^0,v^1)=\mathrm{Span}{(v^5,v^6)}\simeq \C^2.
\end{align*}

Therefore, for some scalar $\lambda_6,\lambda_7,\mu_4,\mu_5\in \C$, we have
\begin{align*}
	&v^5=\lambda_6v^0+\mu_4v^1\\
	&v^6=\lambda_7v^0+\mu_5v^1
\end{align*}
where $\lambda_6\mu_5-\mu_4\lambda_7\neq 0$. As $R_3(f)=R_1(f)$, we have $v^3\in \mathrm{Span}(v^0,v^1,v^2)$ so that (up to relabelling $v^4$ into $v^3$)
\begin{align*}
	F(z)&=\left(1+\lambda_1z+\lambda_2z^2+\lambda_3z^3+\lambda_4z^4+\lambda_5z^5+\lambda_6z^6+\lambda_7(\pi_1z^7+z^8)\right)v^0\\
	&+\left(z^2+\mu_1z^3+\mu_2z^4+\mu_3z^5+\mu_4z^6+\mu_5(\pi_1z^7+z^8)\right)v^1\\
	&+\left(z^3+\nu_1z^4+\nu_2z^5\right)v^2+z^5v^3.
\end{align*}
 In particular, we find that
\begin{align}\label{eq4}
	f_2=&\Big[(2+\cdots)v^0\wedge v^1+(3z+\cdots)v^0\wedge  v^2+(5z^3+\cdots)v^0\wedge v^3+(z^3+\cdots )v^1\wedge v^2+(3z^5+\cdots)v^1\wedge v^3\nonumber\\
	&+z^6(2+\nu_1z)v^2\wedge v^3\Big]=[F_2(z)]
\end{align}
Recall now that $f$ has a branch point of order $k\geq 1$ at $p\in \P^1\setminus\ens{0,\infty}$ if and only if $F_2$ given in \eqref{eq4} has a zero of order $k$ at $p$.
As $P(z)=z^6(2+\nu_1z)$ has exactly one zero with multiplicity $1$ outside of zero, we have a contradiction, as this polynomial $P$ must have three distinct zeroes in $\C\setminus\ens{0}$.
\end{proof}

\begin{cor}
	There does not exist a complete minimal surface in $\R^3$ with exactly $9$ embedded planar ends (and no other ends).
\end{cor}

\section{Some partial results for degree $10$ unbranched null immersions}

Let $f:\P^1\rightarrow \P^3$ a totally ramified contact curve of degree $2d-1\geq 5$. Then $r_1(f)=\deg(f)-3=2d-4$ and for all $p\in  \P^1$, we have
\begin{align*}
    r_1(f,p)\leq \frac{(2d-1)-3}{2}=d-2
\end{align*}
so $f$ has at least two branch points, and $R_1(f,p)=n\,\cdot p+(2d-4-n)\,\cdot q$ for some distinct $p,q\in  \P^1$ and $1\leq n<2d-4$ implies that $n=d-2$, and as we saw earlier, this implies that $f=f_{d}$. Now assume that $f$ has at least three branch points and that $d=5$, so that $\deg(f)=9$ and $r_1(f)=6$. Then $r_1(f,p)\leq d-2=3$ for all $p\in \P^1$, and if $r_1(f,p)=3$ for some $p\in \P^1$, we have $2 r_1(f,p)+3=9=\deg(f)$ so we cannot apply Lemma \ref{abel}. However, if $r_1(f,p)=2$ for some $p\in \P^1$, then $2r_1(f,p)+3=7<9=\deg(f)$ so by Lemma \ref{abel}, we obtain for all $q\in \P^1\setminus\ens{p}$ the inequality
\begin{align*}
	2r_1(f,p)+r_1(f,q)\leq \deg(f)-4=5
\end{align*}
or 
\begin{align*}
	r_1(f,q)\leq 5-4=1
\end{align*}
so $r_1(f,p)=2$ for some $p\in \P^1$ implies that all other branch points have multiplicity $1$, or 
\begin{align*}
	R_1(f)=2\,p+q+r+s+t
 \end{align*}
for some distinct points $q,r,s,t\in \P^1$. Therefore, the remaining admissible divisors are (for some distinct points $p,q,r,s,t,u\in \P^1$)
\begin{align*}
	R_1(f)&=3\,p+q+r+s\\
	&=2\,p+q+r+s+t\\
	&=p+q+r+s+t+u.
\end{align*}
\textbf{Case 1: $R_1(f)=3\,p+q+r+s$.} Then taking $p=0$, and $q=\infty$, we have with the previous notations
\begin{align*}
	f=\left[\left(1+\lambda_1z+\lambda_2z^2+\lambda_3z^3\right)v^0+z^4v^4+z^5v^5+z^6v^6+z^7v^7+\left(\pi_1z^8+z^9\right)v^9\right].
\end{align*}
By the previous arguments, $(v^0,v^4,v^5,v^9)$ and $(v^9,v^7,v^6,v^4)$ are basis of $\C^4$, and there exists $\lambda_j,\mu_j,\nu_j\in \C$ such that up to renaming of vectors, we have
\begin{align*}
	f=&\Big[\left(1+\lambda_1z+\lambda_2z^2+\lambda_3z^3+\lambda_4z^4+\lambda_5z^5+\lambda_9(\pi_1z^8+z^9)\right)v^0+\left(z^4+\mu_5z^5+\mu_9\left(\pi_1z^8+z^9\right)\right)v^4\\
	&+\left(z^5+\nu_9\left(\pi_1z^8+z^9\right)\right)v^5+z^6v^6+z^7v^7+\left(\pi_1z^8+z^9\right)v^9\Big],
\end{align*}
where $(v^0,v^4,v^5,v^9)$ is an orthogonal basis of $\C^4$. Also, notice as $R_3(f)=R_1(f)$ that 
\begin{align*}
   v^6,v^7\in \mathrm{Span}(v^0,v^4,v^5)
\end{align*}
and that $\pi_1=0$ (otherwise, we would have $v^9\in \mathrm{Span}(v^0,v^4,v^5)$ which would imply that $f$ is degenerate, a contradiction by Proposition \ref{nonlinear}). Furthermore, as $(v^0,v^4,v^5,v^9)$ is an orthogonal basis of $\C^4$ and $(v^4,v^6,v^7,v^9)$ is a basis of $\C^4$, this is now manifest that 
\begin{align*}
	\mathrm{Span}(v^0,v^5)=\mathrm{Span}(v^6,v^7)\simeq \C^2
\end{align*}
so we can write
\begin{align*}
	f=&\Big[\left(1+\lambda_1z+\lambda_2z^2+\lambda_3z^3+\lambda_4z^4+\lambda_5z^5+\lambda_6z^6+\lambda_7z^7+\lambda_9z^9\right)v^0\\
	&+\left(z^4+\mu_5z^5+\mu_9z^9\right)v^4+\left(z^5+\nu_6z^6+\nu_7z^7+\nu_9z^9\right)v^5+z^9v^9\Big].
\end{align*}
Then we compute
\begin{align*}
	f_2=&\Big[\left(4+\cdots\right)v^0\wedge v^4+\left(5z+\cdots\right)v^0\wedge v^4+\left(9z^5+\cdots\right)v^0\wedge v^9+\left(z^5+\cdots\right)v^4\wedge v^5\\
	&+z^9\left(5+4\mu_5z\right)v^4\wedge v^9+\left(4z^{10}+\cdots\right)v^5\wedge v^9\Big].
\end{align*}
However, we see that $z^9(5+4\mu_5z)$ must have $2$ zeroes with multiplicity $1$ in $\C\setminus\ens{0}$, but as this polynomial has at most $1$ zero, we have a contradiction.

\textbf{Case 2: $R_1(f)=2\,p+q+r+s+t$.} Taking $p=0$ and $q=\infty$, we obtain
\begin{align*}
	f=\Big[\left(1+\lambda_1z+\lambda_2z^2\right)v^0+z^3v^3+z^4v^4+z^5v^5+z^6v^6+z^7v^7+\left(\pi_1z^8+z^9\right)v^9\Big].
\end{align*}
Making a linear change of variable, we can assume for more notational convenience that $\pi_1=0$. By the same argument and the orthogonalisation process as in case $2$, and as $(v^0,v^3,v^4,v^7)$ and $(v^9,v^7,v^6,v^4)$ are two basis of $\C^4$, we can write $f$ as
\begin{align*}
	f=&\Big[\left(1+\lambda_1z+\lambda_2z^2+\lambda_3z^3+\lambda_4z^4+\lambda_5z^5+\lambda_6z^6+\lambda_7z^7+\lambda_9z^9\right)v^0\\
	&+\left(z^3+\mu_4z^4+\mu_5z^5+\mu_6z^6+\mu_7z^7+\mu_9z^9\right)v^3+\left(z^4+\nu_5z^5+\nu_7z^7\right)v^4+z^7v^7\Big]
\end{align*}
and
\begin{align*}
	f_2=&\Big[\left(3+\cdots\right)v^0\wedge v^3+\left(4z+\cdots\right)v^0\wedge v^4+\left(7z^4+\cdots\right)v^0\wedge v^7\\
	&+\left(z^4+\cdots\right)v^3\wedge v^4+\left(4z^7+\cdots\right)v^3\wedge v^7+z^8\left(3+2\nu_5z\right)v^4\wedge v^7\Big]
\end{align*}
As the polynomial $z^8\left(3+2\nu_5z\right)$ must have $3$ distinct zeroes on $\C\setminus\ens{0}$, while it has at most $1$ zero on $\C\setminus\ens{0}$, we also have a contradiction.

\textbf{Case 3: $R_1(f)=p+q+r+s+t+u$.} Here direct computations seem to become too difficult, although a computer-assisted proof might be possible.

\section{Some partial results for degree $11$ unbranched null immersions}

Let $f:\P^1\rightarrow \P^3$ the corresponding totally ramified contact curve of degree $2d=10$. Then we have $r_1(f)=2d-3=7$, and for all $p\in \P^1$, 
\begin{align*}
	r_1(f,p)\leq d-3=2.
\end{align*}
Therefore, the admissible divisors of $f$ are for some distinct $p,q,r,s,t,u,v\in \P^1$
\begin{align}\label{tooku}
	R_1(f)&=2\,p+2\,q+2\,r+s\nonumber\\
	&=2\,p+2\,q+r+s+t\nonumber\\
	&=2\,p+q+r+s+t+u\nonumber\\
	&=p+q+r+s+t+u+v.
\end{align}
A similar trick as in the proof of Theorem \ref{deg9} permits to rule out the first two divisors. Indeed, if $f$ has two branch points of order $2$ that we take at $z=0$ and $z=\infty$, we obtain with obvious notations
\begin{align*}
	f=\left[\left(1+\lambda_1z+\lambda_2z^2\right)v^0+z^3v^3+z^4v^4+z^5v^5+z^6v^6+z^7v^7+\left(\pi_1z^8+\pi_2z^9+z^{10}\right)v^{10}\right].
\end{align*}
As $R_2(f)=0$ and $R_3(f)=R_1(f)$, we see that both $(v^0,v^3,v^4,v^7)$ and $(v^3,v^6,v^7,v^{10})$ ar basis of $\C^4$. 
Using the same trick as in the proof of Theorem \eqref{deg8} and $R_3(f)=R_1(f)$, we let $\lambda_j,\mu_j,\nu_j\in \C$ be such that $(v^0,v^2,v^3,v^7)$ be an orthogonal basis of $\C^4$ and
\begin{align*}
	f=&\Big[\left(1+\lambda_1z+\lambda_2z^2+\lambda_3z^3+\lambda_4z^4+\lambda_7z^7\right)v^0+\left(z^3+\mu_4z^4+\mu_7z^7\right)v^3+\left(z^4+\nu_7z^7\right)v^4\\
	&+z^5v^5+z^6v^6+z^7v^7+\left(\pi_1z^8+\pi_2z^9+z^{10}\right)v^{10}\Big].
\end{align*}
Recalling that $(v^3,v^6,v^7,v^{10})$ is a basis of $\C^4$ and $(v^0,v^2,v^3,v^7)$ is an \emph{orthogonal} basis of $\C^4$, we deduce that there exists that $v^6,v^{10}\in \mathrm{Span}(v^0,v^4)$. Therefore, there exists $\lambda_6,\lambda_8,\nu_6,\nu_8\in \C$ such that
\begin{align*}
	f=&\Big[\left(1+\lambda_1z+\lambda_2z^2+\lambda_3z^3+\lambda_4z^4+\lambda_6z^6+\lambda_7z^7+\lambda_8\left(\pi_1z^8+\pi_2z^9+z^{10}\right)\right)v^0+\left(z^3+\mu_4z^4+\mu_7z^7\right)v^3\\
	&+\left(z^4+\nu_6z^6+\nu_8\left(\pi_1z^8+\pi_2z^9+z^{10}\right)\right)v^4
	+z^5v^5+z^7v^7\Big]
\end{align*}
Furthermore, we have $v^5\in \mathrm{Span}(v^0,v^3,v^4)$ as $r_3(f,0)=r_1(f,0)=2$, so we finally obtain for some additional $\lambda_5,\mu_5,\nu_5\in \C$
\begin{align*}
	f=&\Big[\left(1+\lambda_1z+\lambda_2z^2+\lambda_3z^3+\lambda_4z^4+\lambda_5z^5+\lambda_6z^6+\lambda_7z^7+\lambda_8\left(\pi_1z^8+\pi_2z^9+z^{10}\right)\right)v^0\\
	&+\left(z^3+\mu_4z^4+\mu_5z^5+\mu_7z^7\right)v^3
	+\left(z^4+\nu_5z^5+\nu_8\left(\pi_1z^8+\pi_2z^9+z^{10}\right)\right)v^4+z^7v^7\Big].
\end{align*}
Now, we compute
\begin{align*}
	f_2&=\Big[\left(3+\cdots\right)v^0\wedge v^3+\left(4z+\cdots\right)v^0\wedge v^4+\left(7z^4+\cdots\right)v^0\wedge v^7+\left(z^4+\cdots\right)v^3\wedge v^4\\
	&\;\,+z^7\left(4+3\mu_4z+2\mu_5z^2\right)v^3\wedge v^7+\left(3z^8+\cdots\right)v^4\wedge v^7\Big]\\
	&=[F_2(z)]
\end{align*}
Now, we see that $F_2$ must have $3$ zeroes with multiplicity in $\C\setminus\ens{0}$ thanks of \eqref{tooku}. However, the coefficient in $v^3\wedge v^7$ is 
\begin{align*}
	z^7\left(4+3\mu_4z+2\mu_5z^2\right)
\end{align*}
which has as most $2$ zeroes in $\C\setminus\ens{0}$, a contradiction.

Therefore, the only possible branch divisors for $f$ are
\begin{align*}
	R_1(f)&=2\,p+q+r+s+t+u\nonumber\\
	&=p+q+r+s+t+u+v
\end{align*}
out of the $15$ possibilities initially.

\section{Remarks on totally ramified contact curves}\label{excluded2}

We saw that when two branch points of a totally ramified contact curve satisfy some algebraic property, then the curve cannot exist. Outside of Lemmas \ref{cont} and \ref{abel}, they correspond to \emph{a subset} of the following situation : there exists distinct $p,q\in \P^1$ such that $a=r_1(f,p)$ and $b=r_1(f,q)$ satisfy
\begin{align}\label{c}
	\mathrm{Card}\left(\ens{0,a+1,a+2,2a+3}\cap \ens{d-(2b+3),d-(b+2),d-(b+1),d}\right)\geq 2,
\end{align}
where we denoted $d=\deg(f)$. Let us consider case by case when this relation does occur. Recall that $2r_1(f,p)+3\leq d$ for all $p\in \P^1$. We summarise the results in the following proposition.
\begin{prop}
	Let $f:\P^1\rightarrow \P^3$ be a totally ramified contact curve of degree $d$ and branch divisor $R_1(f)$ such that \eqref{c} holds for some distinct $p,q\in \P^1$. Then $d$ is \emph{odd} and $f$ is (up to projective equivalence) the element $f_{d'}$ of the family \eqref{fd}, where $d=2d'-1$.
\end{prop}
\begin{proof}
	As there are many cases to treat, we will adopt the notation $(i\, k)\; (j\, l)$ whenever $1\leq i<j\leq 3$ and $1\leq k<l\leq 3$ to say that the $i$-th element (resp $j$-th element) of $\ens{0,a+1,a+2,2a+3}$ corresponds to the $k$-th element (resp. $l$-th element) of $\ens{d-(2b+3),d-(b+2),d-(b+1),d}$. For example, $(1\, 1)\; (2\, 2)$ corresponds to $0=d-(2b+3)$ and $a+1=d-(b+2)$.  This notation will ensure that indeed all cases are included in the forthcoming discussion. For the sake of readability, we write the two collections of indices in \eqref{c} as
	\begin{align}\label{index}
		\begin{pmatrix}
		(1)\\
		(2)\\
		(3)\\
		(4)
		\end{pmatrix}\begin{pmatrix}
		0\\
		a+1\\
		a+2\\
		2a+3
		\end{pmatrix},\quad \begin{pmatrix}
		(1)\\
		(2)\\
		(3)\\
		(4)
		\end{pmatrix}\begin{pmatrix}
		d-(2b+3)\\
		d-(b+2)\\
		d-(b+1)\\
		d
		\end{pmatrix}
	\end{align}

	\textbf{Case 1: $(1\, 1)\; (\ast\,\ast)$}. Then $0=d-(2b+3)$, so $d$ is odd an we replace $d$ by $2d+1$, and $b=d-1$. This implies that
	\begin{align}\label{w1}
		\deg(f)=2d+1,
	\end{align}
	and \eqref{index} becomes
	\begin{align}\label{index11}
		\begin{pmatrix}
		(1)\\
		(2)\\
		(3)\\
		(4)
		\end{pmatrix}\begin{pmatrix}
		0\\
		a+1\\
		a+2\\
		2a+3
		\end{pmatrix},\quad \begin{pmatrix}
		(1)\\
		(2)\\
		(3)\\
		(4)
		\end{pmatrix}\begin{pmatrix}
		0\\
	    d\\
		d+1\\
		2d+1
		\end{pmatrix}.
	\end{align}
	Notice that by Lemma \ref{cont} we have $r_1(f,p)\leq d-1$ for all $p\in \P^1$.
	
	\textbf{Sub-case 1: $(1\, 1)\; (2\, 2)$}. Then $a=b=d-1$, so $f=f_{d+1}$ is the curve given in \eqref{fd}. 
	
	\textbf{Sub-case 2: $(1\, 1)\; (2\, 3)$}. Then $a=d>d-1$, contradiction by Lemma \ref{cont}.
	
	\textbf{Sub-case 3: $(1\, 1)\; (2\, 4)$}. Then $a+1=2d+1$, so $a=2d>\frac{\deg(f)-3}{2}=d-1$, contradiction.
	
	\textbf{Sub-case 4: $(1\, 1)\; (3\, 2)$}. Then $a=d-2$, and as $2(d-2)+3=2d-1<\deg(f)=2d+1$, by have by Lemma \ref{abel}
	\begin{align*}
	2(d-2)+(d-1)=2a+b\leq \deg(d)-4=2d-3,
	\end{align*}
	or $d\leq 2$. Therefore, $\deg(f)=2d+1\leq 5$. However, there are no totally ramified contact curves of degree $3$ and $5$ and more than $2$ branch points by Theorem \ref{TC} (\cite{bryantcurves}).
	
	\textbf{Sub-case 5: $(1\, 1)\; (3\, 3)$}. Then $a=b=d-1$, so $f=f_{d+1}$ (see \eqref{fd}).                          
	
	\textbf{Sub-case 6: $(1\, 1)\; (3\, 4)$}. Then $a=2d-1>d-1$ as $d\geq 1$,  contradiction.    
	
	\textbf{Sub-case 7: $(1\, 1)\; (4\, 2)$}. Then $2a+3=d$, so $d$ is odd, so replacing $d$ by $2d+1$, we obtain
	\begin{align}\label{c1142}
		\deg(f)=2(2d+1)+1=4d+3,\quad a=d-1,\;\, b=2d.
	\end{align}
	Taking $q=0$ and $p=\infty$,  we obtain an expansion
	\begin{align*}
	f=&\Big[\left(1+\lambda_1z+\cdots+\lambda_{2d}z^{2d}\right)v^0+z^{2d+1}v^{2d+1}+z^{2d+2}v^{2d+2}+\cdots+z^{3d+2}v^{3d+2}+z^{3d+3}v^{3d+3}\\
	&+\left(\pi_{3d+4}z^{3d+4}+\cdots+z^{4d+3}\right)v^{4d+3}\Big].
	\end{align*}
	Here, $(v^0,v^{2d+1},v^{2d+2},v^{4d+3})$ and $(v^{2d+2},v^{3d+2},v^{3d+3},v^{4d+3})$ are both basis of $\C^4$. Furthermore, notice that we must have $\pi_j=0$ for all $3d+4\leq j\leq 4d+2$, otherwise $f$ would be degenerate. Now, let $\lambda_j,\mu_j,\nu_j\in \C$ such that $(v^0,v^{2d+1},v^{2d+2},v^{4d+3})$ be an orthogonal basis of $\C^4$ and
	\begin{align*}
	f=&\Big[\left(1+\lambda_1z+\cdots+\lambda_{2d}z^{2d}+\lambda_{2d+1}z^{2d+1}+\lambda_{2d+2}z^{2d+2}+\lambda_{4d+3}z^{4d+3}\right)v^0\\
	&+\left(z^{2d+1}+\mu_{2d+2}z^{2d+2}+\mu_{4d+3}z^{4d+3}\right)v^{2d+1}\\
	&+\left(z^{2d+2}+\nu_{4d+3}z^{4d+3}\right)v^{2d+2}+\cdots+z^{3d+2}v^{3d+2}+z^{3d+3}v^{3d+3}+z^{4d+3}v^{4d+3}\Big]
	\end{align*}
	Now we obtain 
	\begin{align*}
	\mathrm{Span}(v^0,v^{2d+1})=\mathrm{Span}(v^{3d+2},v^{3d+3}),
	\end{align*}
	and as $v^j\in \mathrm{Span}(v^0,v^{2d+1},v^{2d+2})$ for all $2d+3\leq j\leq 3d+3$ we obtain
	\begin{align*}
	f=&\Big[\left(1+\lambda_1z+\cdots+\lambda_{3d+3}z^{3d+3}+\lambda_{4d+3}z^{4d+3}\right)v^0\\
	&+\left(z^{2d+1}+\mu_{2d+2}z^{2d+2}+\cdots+\mu_{3d+3}z^{3d+3}+\mu_{4d+3}z^{4d+3}\right)v^{2d+1}\\
	&+\left(z^{2d+2}+\nu_{2d+3}z^{2d+3}+\cdots+\nu_{3d+1}z^{3d+1}+\nu_{4d+3}z^{4d+3}\right)v^{2d+2}+z^{4d+3}v^{4d+3}\Big].
	\end{align*}
	Finally
	\begin{align*}
	f_2=&\Big[\left((2d+1)+\cdots\right)v^0\wedge v^{2d+1}+\left((2d+2)z+\cdots\right)v^0\wedge v^{2d+2}+\left((4d+3)z^{2d+2}+\cdots\right)v^0\wedge v^{4d+3}\\
	&+\left(z^{2d+2}+\cdots\right)v^{2d+1}\wedge v^{2d+2}+\left((2d+2)z^{4d+3}+\cdots\right)v^{2d+1}\wedge v^{4d+3}\\
	&+z^{4d+4}\left((2d+1)+2d\,\nu_{2d+3}z+\cdots+(d+2)\nu_{3d+1}z^{d-1}\right)v^{2d+2}\wedge v^{4d+3}\Big]=[F_2(z)].
	\end{align*}
	As previously, notice that $r_1(f)-r_1(f,p)-r_1(f,q)=(4d+3)-3-2d-(d-1)=d+1$, so $F_2$ must have $d+1$ roots (with multiplicity) in $\C\setminus\ens{0}$, but $(2d+1)+2d\,\nu_{2d+3}z+\cdots+(d+2)\nu_{3d+1}z^{d-1}$ has degree at most $d-1<d+1$, so we have again a contradiction.
	
	\textbf{Sub-case 8: $(1\,1)\; (4\, 3)$}. Then $2a+3=d+1$, so $d$ is even and we replace $d$ by $2d$, so that
	\begin{align*}
		\deg(f)=4d+1,\quad a=d-1,\quad b=2d-1.
	\end{align*}
	As $a\geq 1$, notice that it implies that $\deg(f)\geq 9$ (also, remark that these multiplicities represent a borderline case of Lemma \ref{abel}). Now take $p=\infty$, $q=0$ and write
	\begin{align*}
	f=&\Big[\left(1+\lambda_1z+\cdots+\lambda_{2d-1}z^{2d-1}\right)v^0+z^{2d}v^{2d}+z^{2d+1}v^{2d+1}+\cdots+z^{3d}v^{3d}+z^{3d+1}v^{3d+1}\\
	&+\left(\pi_{3d+2}z^{3d+2}+\cdots+z^{4d+1}\right)v^{4d+1}\Big].
	\end{align*}
	We first remark as $R_3(f)=R_1(f)$ that $(v^0,v^{2d},v^{2d+1},v^{4d+1})$ is a basis of $\C^4$ so $\pi_j=0$ for all $3d+2\leq j\leq 4d$, and $(v^{2d},v^{3d},v^{3d+1},v^{4d+1})$ is also a basis of $\C^4$. Let $\lambda_j,\mu_j,\nu_j\in \C$ such that
	\begin{align*}
	f=&\Big[\left(1+\lambda_1z+\cdots+\lambda_{2d-1}z^{2d-1}+\lambda_{2d}z^{2d}+\lambda_{2d+1}z^{2d+1}+\lambda_{4d+1}z^{4d+1}\right)v^0\\
	&+\left(z^{2d}+\mu_{2d+1}z^{2d+1}+\mu_{4d+1}z^{4d+1}\right)v^{2d}+\left(z^{2d+1}+\nu_{4d+1}z^{4d+1}\right)v^{2d+1}\\
	&+\cdots
	+z^{3d}v^{3d}+z^{3d+1}v^{3d+1}
	+z^{4d+1}v^{4d+1}\Big]
	\end{align*}
	and $(v^0,v^{2d+1},v^{2d+2},v^{4d+1})$ be an orthogonal basis of $\C^4$. Therefore, by the previous argument we have
	\begin{align*}
	\mathrm{Span}(v^0,v^{2d+1})=\mathrm{Span}(v^{3d},v^{3d+1})\simeq \C^2.
	\end{align*}
	We deduce as $v^j\in \mathrm{Span}(v^0,v^{2d},v^{2d+1})$ for all $2d+2\leq j\leq 3d+1$ that there exists $\lambda_j,\mu_j,\nu_j\in \C$ such that
	\begin{align*}
	f=&\Big[\left(1+\lambda_1z+\cdots+\lambda_{3d+1}z^{3d+1}+\lambda_{4d+1}z^{4d+1})\right)v^0\\
	&+\left(z^{2d}+\mu_{2d+1}z^{2d+1}+\mu_{2d+2}z^{2d+2}+\cdots+\mu_{3d-1}z^{3d-1}+\mu_{4d+1}z^{4d+1}\right)v^{2d}\\
	&+\left(z^{2d+1}+\nu_{2d+2}z^{2d+2}+\cdots+\nu_{3d+1}z^{3d+1}+\nu_{4d+1}z^{4d+1}\right)v^{2d+1}+z^{4d+1}v^{4d+1}\Big].
	\end{align*}
	Now, we compute
	\begin{align*}
	f_2=&\Big[\left(2d+\cdots\right)v^0\wedge v^{2d}+\left((2d+1)z+\cdots\right)v^{0}\wedge v^{2d+1}+\left((4d+1)z^{2d+1}+\cdots\right)v^0\wedge v^{4d+1}\\
	&+\left(z^{2d+1}+\cdots\right)v^{2d}\wedge v^{2d+1}+z^{4d+1}\left((2d+1)+2d\,\mu_{2d+1}z+\cdots+(d+2)\mu_{3d-1}z^{d-1}\right)v^{2d}\wedge v^{4d+1}\\
	&+\left(2d\,z^{4d+2}+\cdots\right)v^{2d+1}\wedge v^{4d+1}\Big]=[F_2(z)].
	\end{align*}
	Now, we have $r_1(f)=\deg(f)-3=4d-2$, so $F_2$ admits exactly $4d-2-(2d-1)-(d-1)=d$ zeroes on $\C\setminus\ens{0}$, while $(2d+1)+2d\,\mu_{2d+1}z+\cdots+(d+2)\mu_{3d}z^{d-1}$ has degree at most $d-1<d$, so we have a contradiction.
	
	\textbf{Sub-case 9: $(1\, 1)\; (4\, 4)$}. Then $a=b=d-1$, so $f=f_{d+1}$ given  by \eqref{fd}.
	
	This concludes the proof of the $(1\,1 )\;(\ast\,\ast)$ case.
	
	\textbf{Case 2: $(1\,2)\;(\ast\,\ast)$ or $(1\,3)\; (\ast\,\ast)$.} In both cases, $b\geq d-2>\frac{d-3}{2}$, contradiction.
	
	\textbf{Case 3: $(2\, 1)\; (\ast\,\ast)$}. Then $a+1=d-(2b+3)$ or 
	\begin{align}\label{w2}
	   a+2b=d-4.
	\end{align}
	
	\textbf{Sub-case 1: $(2\,1)\; (3\, 2)$}. This implies that $a+b=d-4$, so $b=0$ by \eqref{w2}, contradiction.
	
	\textbf{Sub-case 2: $(2\, 1)\; (3 \,3)$}. Then $a+b=d-3<a+2b=d-4$, absurd.
	
	\textbf{Sub-case 3: $(2\, 1)\; (3\, 4)$}. Then $a=d-2>\frac{d-3}{2}$, absurd.
	
	\textbf{Sub-case 4: $(2\, 1)\; (4\, 2)$}. Then $2a+b=d-5$, and
	\begin{align}\label{s2142}
		\left\{\begin{alignedat}{1}
		2a+b=d-5\\
		a+2b=d-4
		\end{alignedat}\right.
	\end{align}
	implies that  $d=0\;(\text{mod}\; 3)$, so replacing $d$ by $3d$, we obtain
	\begin{align*}
		\deg(f)=3d,\quad a=d-2,\quad b=d-1.
	\end{align*}
	Taking $p=\infty$, $q=0$, so that 
	\begin{align*}
	a=r_1(f,p)=d-2,\quad b=r_1(f,q)=d-1,
	\end{align*}
	we can write $f$ as
	\begin{align*}
	f=&\Big[\left(1+\lambda_1z+\cdots+\lambda_{d-1}z^{d-1}\right)v^0+z^dv^d+z^{d+1}v^{d+1}+\cdots+z^{2d+1}v^{2d+1}+z^{2d+2}v^{2d+2}\\
	&+\left(\pi_{2d+3}z^{2d+3}+\cdots+z^{3d}\right)v^{3d}\Big].
	\end{align*}
	Now, notice that $(v^0,v^d,v^{d+1},v^{2d+1})$ and $(v^{d+1},v^{2d},v^{2d+1},v^{3d})$ are two basis of $\C^4$. Using the same method as before, let $\lambda_j,\mu_j,\nu_j\in \C$ such that
	\begin{align*}
	f=&\Big[\left(1+\lambda_1z+\cdots\lambda_{d-1}z^{d-1}+\lambda_dz^d+\lambda_{d+1}z^{d+1}+\lambda_{2d+1}z^{2d+1}\right)v^0+\left(z^d+\mu_{d+1}z^{d+1}+\mu_{2d+1}z^{2d+1}\right)v^d\\
	&+\left(z^{d+1}+\nu_{2d+1}z^{2d+1}\right)v^{d+1}+\cdots+z^{2d}v^{2d}+z^{2d+1}v^{2d+1}+\left(\pi_{2d+2}z^{2d+2}+\cdots+z^{3d}\right)v^{3d}\Big].
	\end{align*}
	and such that $(v^0,v^d,v^{d+1},v^{2d+1})$ be an orthogonal basis of $\C^4$. As $(v^{d+1},v^{2d},v^{2d+1},v^{3d})$ is also an orthogonal basis of $\C^4$, we obtain
	\begin{align*}
	\mathrm{Span}(v^0,v^d)=\mathrm{Span}(v^{2d},v^{3d})\simeq \C^2.
	\end{align*}
	Therefore, as $v^j\in \mathrm{Span}(v^0,v^{d},v^{d+1})$ for all $d+2\leq j\leq 2d$ we have
	\begin{align*}
	f=&\Big[\left(1+\lambda_1z+\cdots+\lambda_{2d+1}z^{2d+1}+\lambda_{2d+2}\left(\pi_{2d+2}z^{2d+2}+\cdots+z^{3d}\right)\right)v^0\\
	&+\left(z^d+\mu_{d+1}z^{d+1}+\cdots+\mu_{2d+1}z^{2d+1}+\mu_{2d+2}\left(\pi_{2d+2}z^{2d+2}+\cdots+z^{3d}\right)\right)v^d\\
	&+\left(z^{d+1}+\nu_{d+2}z^{d+2}+\cdots+\nu_{2d-1}z^{2d-1}+\nu_{2d+1}z^{2d+1}\right)v^{d+1}+z^{2d+1}v^{2d+1}\Big].
	\end{align*}
	Finally, we compute
	\begin{align*}
	f_2&=\Big[\left(d+\cdots\right)v^0\wedge v^d+\left((d+1)z+\cdots\right)v^0\wedge v^{d+1}+\left((2d+1)z^{d+1}+\cdots\right)v^0\wedge v^{2d+1}
	+\left(z^{d+1}+\cdots\right)v^d\wedge v^{d+1}\\
	&\quad\;\,+\left((d+1)z^{2+1d}+\cdots\right)v^d\wedge v^{2d+1}+z^{2d+2}\left(d+(d-1)\nu_{d+2}+\cdots+2\nu_{2d-1}z^{d-1}\right)v^{d+1}\wedge v^{2d+1}\Big]\\
	&=[F_2(z)]
	\end{align*}
	Now, we have $r_1(f,p)+r_1(f,q)=2d-3$, while $r_1(f)=\deg(f)-3=3d-3$, so $F_2(z)$ admits exactly $3d-3-(2d-3)=d$ zeroes with multiplicity on $\C\setminus\ens{0}$. However, $P(z)=d+(d-1)\nu_{d+2}+\cdots+2\nu_{2d-1}z^{d-1}$ has degree at most $d-1$, so we have a contradiction.
	
	\textbf{Sub-case 5: $(2\,1)\; (4\,3)$}. Then
	\begin{align*}
		\left\{\begin{alignedat}{1}
		a+2b=d-4\\
		2a+b=d-4
		\end{alignedat}\right.
	\end{align*}
	so $d=1\;(\text{mod}\; 3)$. Therefore, replace $d$ by $3d+1$, we find
	\begin{align*}
		\deg(f)=3d+1,\quad a=b=d-1.
	\end{align*}
	Taking $p=0$ and $q=\infty$ we obtain the expansion
	\begin{align*}
	f=&\Big[\left(1+\lambda_1z+\cdots+\lambda_{d-1}z^{d-1}\right)v^0+z^dv^d+z^{d+1}v^{d+1}+\cdots+z^{2d}v^{2d}+z^{2d+1}v^{2d+1}\\
	&+\left(\pi_{2d+2}z^{2d+2}+\cdots+z^{3d+1}\right)v^{3d+1}\Big]=[F(z)].
	\end{align*}
	Here, notice that $(v^0,v^d,v^{d+1},v^{2d+1})$ and $(v^d,v^{2d},v^{2d+1},v^{3d+1})$ are basis of $\C^4$. Now, making a change of basis such that $(v^0,v^d,v^{d+1},v^{2d+1})$ becomes orthogonal, we can write $f$ as
	\begin{align*}
	f=&\Big[\left(1+\lambda_1z+\lambda_{d-1}z^{d-1}+\lambda_dz^{d}+\lambda_{d+1}z^{d+1}+\lambda_{2d+1}z^{2d+1}\right)v^0+\left(z^d+\mu_{d+1}z^{d+1}+\mu_{2d+1}z^{2d+1}\right)v^d\\
	&+\left(z^{d+1}+\nu_{2d+1}z^{2d+1}\right)v^{d+1}+\cdots+z^{2d}v^{2d}+z^{2d+1}v^{2d+1}+\left(\pi_{2d+2}z^{2d+2}+\cdots+z^{3d+1}\right)v^{3d+1}\Big].
	\end{align*}
	Recalling that $(v^d,v^{2d},v^{2d+1},v^{3d+1})$ is basis of $\C^4$ and as $(v^0,v^d,v^{d+1},v^{2d+1})$ is an orthogonal basis of $\C^4$, while $v^{j}\in \mathrm{Span}(v^0,v^d,v^{d+1})$ for all $d+2\leq j\leq 2d$, we have
	\begin{align*}
	\mathrm{Span}(v^0,v^{d+1})=\mathrm{Span}(v^{2d},v^{3d+1}),
	\end{align*}
	so we obtain the expansion
	\begin{align}\label{mu2d}
	f=&\Big[\left(1+\lambda_1z+\cdots+\lambda_{2d+1}z^{2d+1}+\lambda_{2d+2}\left(\pi_{2d+2}z^{2d+2}+\cdots+z^{3d+1}\right)\right)v^0\nonumber\\
	&+\left(z^d+\mu_{d+1}z^{d+1}+\cdots+\mu_{2d-1}z^{2d-1}+\mu_{2d+1}z^{2d+1}\right)v^d\nonumber\\
	&+\left(z^{d+1}+\nu_{d+2}z^{d+2}+\cdots+\nu_{2d}z^{2d}+\nu_{2d+1}z^{2d+1}+\nu_{2d+2}\left(\pi_{2d+2}z^{2d+2}+\cdots+z^{3d+1}\right)\right)v^{d+1}+z^{2d+1}v^{2d+1}\Big]\nonumber\\
	&=[F(z)].
	\end{align}
	Now, notice that $r_1(f)=\mathrm{deg}(f)-3=3d-2$, and that
	\begin{align*}
	r_1(f)-r_1(f,p)-r_1(f,q)=r_1(f)-a-b=3d-2-2(d-1)=d,
	\end{align*}
	so the function (up to renaming $v^1=v^d$, $v^2=v^{d+1}$ and $v^{3}=v^{2d+1}$)
	\begin{align*}
	F_2(z)=F(z)\wedge F'(z)=\sum_{0\leq i<j\leq 3}^{}F_{i,j}(z)v^i\wedge v^j
	\end{align*}
	must be such that $F_{i,j}$ admits exactly $d$ zeroes with multiplicity in $\C\setminus\ens{0}$.
	However, we compute (notice that $\mu_{2d}=0$ in \eqref{mu2d})
	\begin{align*}
	z^{1-d}F_2(z)&=\left(d+\cdots\right)v^0\wedge v^d+\left((d+1)z+\cdots\right)v^0\wedge v^{d+1}+\left((2d+1)z^{d+1}+\cdots\right)v^0\wedge v^{2d+1}\\
	&+\left(z^{d+1}+\cdots\right)v^d\wedge v^{d+1}+z^{2d+1}\left((d+1)+d\,\mu_{d+1}z+\cdots+2\mu_{2d-1}z^{d-1}\right)v^d\wedge v^{2d+1}\\
	&+\left(z^{2d+2}+\cdots\right)v^{d+1}\wedge v^{2d+1}
	\end{align*}
	but $z^{2d+1}\left((d+1)+d\mu_{d+1}z+\cdots+2\mu_{2d-1}z^{d-1}\right)$ admits at most $d-1$ zeroes (with multiplicity) on $\C\setminus\ens{0}$, while it must admit exactly $d$ zeroes (with multiplicity), so we have a contradiction.
	
	\textbf{Sub-case 6: $(2\,1)\; (4\,4)$}. Then 
	\begin{align*}
		\left\{\begin{alignedat}{1}
		a+2b=d-4\\
		2a+3=d
		\end{alignedat}\right.
	\end{align*}
	so
	\begin{align*}
		b=\frac{1}{2}\left(d-4-\frac{d-3}{2}\right)=\frac{d-5}{4}
	\end{align*}
	which implies that $d=1\; (\text{mod}\;4)$. Therefore, replacing $d$ by $4d+1$, we obtain
	\begin{align*}
		\deg(f)=4d+1,\quad a=2d-1,\quad b=d-1
	\end{align*}
	which corresponds to the impossible case \textbf{Sub-case 8: $(1\, 1)\; (4\, 3)$}, where $a$ and $b$ are exchanged.
	
	\textbf{Case 4: $(2\,2)\;(\ast\,\ast)$}. Here $a+b=d-3$, and $a,b\leq \frac{d-3}{2}$ implies that $a=b=\frac{d-3}{2}$, corresponding to the cases $(2\,2)\; (3\, 3)$ and $(2\, 2)\; (4\,4)$, while the other cases are excluded by the inequality of Lemma \ref{cont}. Therefore, $(2\,2)\; (3\,4)$ and $(2\,2)\; (4\,3)$ are impossible.
	
	\textbf{Case 5: $(2\,3)\; (\ast\,\ast)$}. Here we have
	\begin{align*}
		a+b=d-2>2\left(\frac{d-3}{2}\right),
	\end{align*}
	a contradiction by Lemma \ref{cont}.
	
	\textbf{Case 6: $(3\,1)\; (\ast\,\ast)$}. Then $a+2b=d-5$.
	
	\textbf{Sub-case 1: $(3\,1)\; (4\,2)$}. This is equivalent to
	\begin{align}
		\left\{\begin{alignedat}{1}
		&a+2b=d-5\\
		&2a+b=d-5
		\end{alignedat}\right.
	\end{align}
	so $a=b=\dfrac{d-5}{3}$, and $d=2\;(\text{mod}\;3)$. Replacing $d$ by $3d+2$, we have
	\begin{align*}
		\deg(f)=3d+2,\quad a=b=d-1.
	\end{align*}
	Taking $p=0$ and $q=\infty$, we have an expansion
	\begin{align*}
	f=&\Big[\left(1+\lambda_1z+\cdots+\lambda_{d-1}z^{d-1}\right)v^0+z^dv^d+z^{d+1}v^{d+1}+\cdots+z^{2d+1}v^{2d+1}+z^{2d+2}v^{2d+2}\\
	&+\left(\pi_{2d+3}z^{2d+3}+\cdots+z^{3d+2}\right)v^{3d+2}\Big].
	\end{align*}
	Here, $(v^0,v^d,v^{d+1},v^{2d+1})$ and $(v^{d+1},v^{2d+1},v^{2d+2},v^{3d+2})$ are basis of $\C^4$. Now, let $\lambda_j,\mu_j,\nu_j\in \C$ be such that
	\begin{align*}
	f=&\Big[\left(1+\lambda_1z+\cdots+\lambda_{d-1}z^{d-1}+\lambda_dz^d+\lambda_{d+1}z^{d+1}+\lambda_{2d+1}z^{2d+1}\right)v^0+\left(z^d+\mu_{d+1}z^{d+1}+\mu_{2d+1}z^{2d+1}\right)v^d\\
	&+\left(z^{d+1}+\nu_{2d+1}z^{2d+1}\right)v^{d+1}+\cdots+z^{2d+1}v^{2d+1}+z^{2d+2}v^{2d+2}
	+\left(\pi_{2d+3}z^{2d+3}+\cdots+z^{3d+2}\right)v^{3d+2}\Big].
	\end{align*}
	and $(v^0,v^{d},v^{d+1},v^{2d+1})$ be an orthogonal basis of $\C^4$. As $(v^{d+1},v^{2d+1},v^{2d+2},v^{3d+2})$ is also a basis of $\C^4$ this is now manifest that
	\begin{align*}
	\mathrm{Span}(v^0,v^d)=\mathrm{Span}(v^{2d+2},v^{3d+2})\simeq \C^2.
	\end{align*}
	Using also the relation $v^j\in \mathrm{Span}(v^0,v^{d},v^{d+1})$ for all $d+2\leq j\leq 2d$ coming from $R_3(f)=R_1(f)$, we obtain
	\begin{align*}
	f=&\Big[\left(1+\lambda_1z+\cdots+\lambda_{2d+2}z^{2d+2}+\lambda_{2d+3}\left(\pi_{2d+3}z^{2d+3}+\cdots+z^{3d+2}\right)\right)v^0\\
	&+\left(z^d+\mu_{d+1}z^{d+1}+\cdots+\mu_{2d+2}z^{2d+2}+\mu_{2d+3}\left(\pi_{2d+3}z^{2d+3}+\cdots+z^{3d+2}\right)\right)v^d\\
	&+\left(z^{d+1}+\nu_{d+2}z^{d+2}+\cdots+\nu_{2d}z^{2d}+\nu_{2d+1}z^{2d+1}\right)v^{d+1}+z^{2d+1}v^{2d+1}\Big].
	\end{align*}
	Therefore, we finally compute
	\begin{align*}
	f_2=&\Big[\left(d+\cdots\right)v^0\wedge v^d+\left((d+1)z+\cdots\right)v^0\wedge v^{d+1}+\left((2d+1)z^{d+1}+\cdots\right)v^0\wedge v^{2d+2}\\
	&+\left(z^{d+1}+\cdots\right)v^d\wedge v^{d+1}+\left((d+2)z^{2d+1}+\cdots\right)v^d\wedge v^{2d+1}\\
	&+z^{2d+2}\left(d+(d-1)\nu_{d+2}z+\cdots+\nu_{2d}z^{d-1}\right)v^{d+1}\wedge v^{2d+1}\Big]=[F_2(z)].
	\end{align*}
	Now, $f$ has exactly $r_1(f)=\deg(f)-3=3d-1$ branch points (with multiplicity), while $r_1(f,p)+r_1(f,q)=2d-2$, so $F_2$ admits exactly $3d-1-(2d-2)=d+1$ zeroes on $\C\setminus\ens{0}$ while $d+(d-1)\nu_{d+2}z+\cdots+\nu_{2d}z^{d-1}$ has degree at most $d-1<d+1$ so we have a contradiction.
	
	\textbf{Sub-case 2: $(3\,1)\; (4\,3)$}. This yields the system
	\begin{align*}
		\left\{\begin{alignedat}{1}
		a+2b=d-5\\
		2a+b=d-4
		\end{alignedat}\right.
	\end{align*}
	which corresponds to the system \eqref{s2142}, where $a$ and $b$ are exchanged, so this case is impossible (it corresponds to the case where $d=0\; (\text{mod}\; 3)$).
	
	\textbf{Sub-case 3: $(3\,1)\; (4\,4)$}. As $2a+3=d$, the curve $f$ has odd degree $d$, and
	\begin{align*}
		b=\frac{1}{2}\left(d-5-\frac{d-3}{2}\right)=\frac{d-7}{4},
	\end{align*}
	so $d=3\; (\text{mod}\;4)$, and replacing $d$ by $4d+3$, we obtain
	\begin{align*}
		\deg(f)=4d+3,\quad a=2d, \quad b=d-1
	\end{align*}
	which corresponds to \textbf{Case 1, Sub-case 7} $(1\, 1)\; (4\,2)$ where $a$ and $b$ are exchanged (see \eqref{c1142}) so this case is also excluded.
	
	\textbf{Case 7: $(3\,2)\; (\ast\,\ast)$}. Here we have $a+b=d-4$.
		
	\textbf{Sub-case 1: $(3\,2)\;(4\,3)$}. Then $2a+b=d-4$ and $a+b=d-4$, which is absurd as $a\geq 1$.
	
	\textbf{Sub-case 2: $(3\, 2)\; (4\, 4)$}. Then $d$ is odd, so replace $d$ by $2d+1$ to obtain $a=d-1$ and $b=d-2$. However, by the argument of \textbf{Sub-case 4: $(1\, 1)\; (3\, 2)$}., this implies that $\deg(f)=2d+1\leq 5$, a contradiction.
	
	\textbf{Case 8: $(3\,3)\;(4\,4)$.} Then $d=2d'-1$ is odd, $a=b=d'-2$ and $f=f_{d'}$ given by \eqref{fd}.
	
	This completes the proof of the proposition.
\end{proof}


\section{On the examples in the literature}\label{literature}

In \cite{peng0}, it is claimed that there exists minimal surfaces with any odd-number $2d+1\geq 9$ of embedded flat ends. However, the paper was never published to my knowledge, and the papers which actually appeared were the following ones: \cite{peng1} and \cite{peng2}. There is a family of examples given for even and odd number of ends, but they fail to have the asserted properties. We treat the case with an odd number of ends, as we already know examples of minimal surfaces with an even number (necessarily larger than $4$, see \cite{kusnerpacific}, \cite{bryant3}, \cite{bryantcurves}) of embedded flat ends. 

Now, let $\Sigma$ be a closed Riemann surface, $p_1,\cdots,p_n\in \Sigma$ be $n\geq 1$ distinct points, and $\phi:\Sigma\setminus\ens{p_1,\cdots,p_n}\rightarrow\R^3$ be a complete minimal surface of finite total curvature and $g:\Sigma\setminus\ens{p_1,\cdots,p_n}\rightarrow\P^1$ its the Gauss map. It is a classical fact (\cite{dierkes}) that $g$ extends continuously at branched points $p_1,\cdots,p_n$. As $\phi$ is minimal, $g$ is a harmonic map so it extends analytically on $\Sigma$. Then this is easy to see that the total curvature of $\phi$ is given as below
\begin{align*}
	C(\phi)=\int_{\Sigma}K_hd\mathrm{vol}_{h}=-4\pi \deg(g),
\end{align*} 
where $h$ is the induced metric of $\phi$ on $\Sigma\setminus\ens{p_1,\cdots,p_n}$.

	There is a well-defined notion of \emph{order} of an en of a complete minimal surface near an end. Fix some $1\leq j\leq n$. By the Weierstrass parametrisation and as $\phi$ is \emph{complete}, for every complex chart $(z,U)$ such that $p_j\in U$, there exists an integer $m+1\geq 2$, such that  for some non-zero constant $A_0\in \C^3\setminus\ens{0}$ (verifying $\s{\vec{A}_0}{\vec{A}_0}=0$)
	\begin{align}\label{embedded}
	\p{z}\phi=\frac{\vec{A_0}}{z^{m+1}}+O(|z|^{-m}).
	\end{align}
	The integer $m\geq 1$ is called the multiplicity of the end $p_j$, and does not depend on the chart. We say that $p_j$ is an embedded end if for all chart $(U,z)$ of sufficiently small enough domain, the restriction $\phi|U\setminus\ens{p_j}\rightarrow\R^3$ is an embedding. If the end is an embedded, we have in particular $m=1$ in \eqref{embedded}. Furthermore, if $m=1$, there exists $a\in \R$ such that after rotation
	\begin{align*}
	\phi(z)=\Re\left(\frac{\vec{A}_0}{z}\right)+(0,0,a)\log|z|+O(|z|.)
	\end{align*}
	and we call $a\in \R$ the \emph{logarithmic growth} of the end $p_j$. If $a\neq 0$, we say that the ends $p_j$ is of \emph{catenoid type} and if  $a=0$ we say that the end is \emph{flat} or \emph{planar}. One can check that the inversion of a complete minimal surface in $\R^3$ is a smooth Willmore \emph{immersion} (without branch) points if and only if its ends are embedded and flat. For example, the inversion of the catenoid is not smooth (and not even $C^{1,1}$), as its ends are embedded but not planar.

Now, if $\Sigma$ has genus $\gamma$ and the ends of $\phi$ are embedded, by the Jorge-Meeks formula (\cite{jorge}), we have
\begin{align*}
	\deg(g)=-\frac{1}{2}\chi(\Sigma)+n=\gamma-1+n.
\end{align*}
In particular, if $\Sigma=\P^1$ we obtain
\begin{align*}
	\deg(g)=n-1.
\end{align*}
Now, the Weierstrass parametrisation of the minimal surfaces in \cite{peng1} is given once we identify $g$ and its stereographic projection by
\begin{align}\label{cont1}
	g_{n,m}(z)=\frac{Q_n(z)}{P_{n,m}(z)},\quad 
	\omega_{n,m}=\frac{P_{n,m}(z)^2}{z^2(z^n-1)^2(z^n-\lambda)^2}dz,
\end{align}
where
\begin{align*}
	P_{n,m}(z)=z^m(z^n-\lambda),\quad  Q_n(z)=(z^n-a)(z^n-b)
\end{align*}
for any $n\geq 4$ and $2\leq m\leq n-1$ such that $2m\neq n+1$, and $a,b,c,\lambda\in \C$ are four distinct points different from $0$ and $1$.
Therefore, we have
\begin{align*}
	g_{n,m}(z)=\frac{(z^n-a)(z^n-b)}{z^m(z^n-c)},\quad \omega_{n,m} =\frac{z^{2m-2}(z^n-c)^2}{(z^n-1)^2(z^n-\lambda)^2}dz,
\end{align*}
and it is claimed that 
\begin{align*}
	\vec{\Psi}_{n,m}(z)=\Re\left(\int_{\ast}^{z}(1-g_{n,m}^2)\,\omega_{n,m},i(1+g_{n,m}^2)\,\omega_{n,m},2g_{n,m}\,\omega_{n,m}\right)
\end{align*}
is a complete minimal surfaces with $2n+1$ embedded planar ends and total curvature $-4\pi(2n)$. As the map $g:\P^1\rightarrow\P^1$ has degree $2n$, if $\vec{\Psi}_{n,m}$ is a well-defined defined minimal surface, we obtain
\begin{align*}
	C(\vec{\Psi}_{n,m})=-4\pi \deg(g_{n,m})=-4\pi(2n).
\end{align*}
However, the $1$-form $\omega$ should have zeroes of order $2m$ when $g$ has poles of order $m$, but in $0$, we remark that $\omega$ has only a pole of order $2m-2$ at $0$ (notice the cancellation between $z^2$ and $P_{n,m}^2$), so $\phi$ cannot be a minimal surface with embedded flat ends. 

We now check explicitly in the simplest example of the family with $n=4$ and $m=2$ that $\vec{\Psi}_{4,2}$ fails to be a complete minimal surface with embedded planar ends.

Indeed, fix $P,Q$ two non-zero relatively prime polynomial functions, and let $a_1,\cdots,a_n\in\C$ be $n$ fixed points. Then the Weierstrass data,
\begin{align*}
	g(z)=\frac{Q(z)}{P(z)},\quad \omega=\frac{P(z)^2}{\prod_{j=1}^n(z-a_j)^2}dz
\end{align*}
gives a complete minimal surface with $n$ embedded ends if and only if $(g,\omega)$ solve the period problem and $P(a_j)Q(a_j)\neq 0$ for all $j=1,\cdots,n$. The period problem, which is equivalent to having the associate minimal surface, corresponds to the conditions
\begin{align*}
	\Re\left(\int_{\gamma}\omega\right)=\Re\left(\int_{\gamma}g\,\omega\right)=\Re\left(\int_{\gamma}g^2\,\omega\right)=0
\end{align*}
for all closed curve $\gamma \subset \Sigma\setminus\ens{p_1,\cdots,p_n}$. If we require to have furthermore planar ends, this condition is equivalent to the absence of residues of the $\C^3$-valued $1$-form $(\omega,g\,\omega,g^2\,\omega)$.
We refer to \cite{jorge} for more details about these definitions. 

For example, the example of complete minimal surface with $2n\geq 4$ embedded planar ends given by Kusner in \cite{kusnerpacific} are given by the following data for any $n\geq 2$ as
\begin{align*}
	g_n(z)=\frac{z^{n-1}(z^n-s_n)}{s_nz^n+1},\quad \omega_n=\frac{i(s_nz^n+1)^2}{(z^{2n}+r_nz^n-1)^2}dz
\end{align*}
where $s_n=\sqrt{2n-1}$, and $r_n=\dfrac{2s_n}{n-1}$. We see that indeed $\omega_n$ has a zero of order $2m$ at poles of $g_n$ of order $m$, and the ends of 
\begin{align}\label{ex}
	\phi_n(z)&=\Re\left(\int_{\ast}^{z}(1-g_n^2)\,\omega_n,i(1+g_n^2)\,\omega_n,2g_n\,\omega_n\right)\nonumber\\
	&=\Re\left(\frac{i}{z^{2n}+r_nz^n-1}\left(z^{2n-1}-z,-i(z^{2n-1}+z),\frac{n-1}{n}(z^{2n}+1)\right)\right)
\end{align} 
corresponds to the $2n$ distinct zeroes of $z^{2n}+r_nz^n-1$. Let
\begin{align*}
	R(z)=z^{2n}+r_nz^n-1=\left(z^n+\frac{r_n}{2}\right)^2-\left(1+\frac{r_n^2}{4}\right)
\end{align*}
then
\begin{align*}
	R'(z)=2nz^{n-1}\left(z^n+\frac{r_n}{2}\right)
\end{align*}
so the zeroes of $R'$ are $0$ and the $n$-roots of $-\dfrac{r_n}{2}$. Now, $R(0)=-1\neq 0$ and if $z^n=-\dfrac{r_n}{2}$, we have
\begin{align*}
	R(z)=-\left(1+\frac{r_n^2}{4}\right)<0
\end{align*}
so $R$ and $R'$ have no common zero.

Let $\phi:\P^1\setminus\ens{a_1,\cdots, a_n}\rightarrow\R^3$ a minimal surface with $n$ embedded planar ends $a_1,\cdots,a_n\in\C=\P^1\setminus\ens{\infty}$. Then there exists a meromorphic null immersion $f:\P^1\rightarrow\C^3$ such that $\phi=\Re(f)$ and
\begin{align*}
	f(z)=v^0+\frac{v^1}{z-a_1}+\cdots +\frac{v^n}{z-a_n}
\end{align*}
for some $v^0\in \C^3$, and $v^1\cdots, v^n\in \C^{3}\setminus\ens{0}$. In particular, if $\phi:\P^1\setminus\ens{a_1,\cdots, a_n}\rightarrow\R^3$ is a minimal surface with $n$ embedded planar ends $a_1,\cdots,a_n\in\C$, with $\phi=\Re(f)$, the function $F:\C\rightarrow\C^3$ given by
\begin{align*}
	F(z)=f(z)\prod_{j=1}^{n}(z-a_j)=v^0\prod_{j=1}^{n}(z-a_j)+\sum_{i=1}^{n}v^i\prod_{j=1,j\neq i}^{n}(z-a_j)
\end{align*}
is a polynomial function (with values in $\C^3$) of degree at most $n$ (exactly equal to $n$ if $v^0\neq 0$, and to $n-1$ if $v^0=0$). With the example in \eqref{ex}, we indeed have with obvious notations
\begin{align*}
	F(z)=i\left(z^{2n-1}-z,-i(z^{2n-1}+z),\frac{n-1}{n}(z^{2n}+1)\right)
\end{align*}
which is indeed a polynomial function with values in $\C^3$ of degree $2n$.

Therefore, it suffices us to show that the corresponding $F_{n,m}$ given by $\vec{\Psi}_{n,m}$ does not enjoy this property to prove that the proposed minimal surface is not a minimal surface with $2n+1$ embedded planar ends. Notice that we also proved that the case $n=4$ is impossible, and we will check indeed that the given examples cannot work this specific case.

Taking the simplest example where $n=4$ and $m=2$ (notice that $4=2m\neq n+1=5$) we can find an admissible solution $(a,b,c,\lambda)$ of parameters  given as
\small
\begin{align*}
    \left\{\begin{alignedat}{1}
	&a=\frac{488 \, \sqrt{15} - 3 \, \sqrt{215208 \, \sqrt{15} + 833497} + 1890}{3 \, {\left(8 \, \sqrt{15} + 31\right)}}=-0.502000420331517
	\cdots,\\
	&b=-\frac{9 \, \sqrt{15} \sqrt{215208 \, \sqrt{15} + 833497} - 8983 \, \sqrt{15} + 33 \, \sqrt{215208 \, \sqrt{15} + 833497} - 34791}{3 \, {\left(5 \, \sqrt{15} \sqrt{215208 \, \sqrt{15} + 833497} - 6747 \, \sqrt{15} + 21 \, \sqrt{215208 \, \sqrt{15} + 833497} - 26131\right)}}=41.1579116001055\cdots
	,\\
	&c=-\frac{7 \, \sqrt{15} + 27}{\sqrt{15} + 5}=-6.09838667696593
	\cdots,\\
	&\lambda=-8 \, \sqrt{15} - 31=-61.9838667696593
	\cdots.
	\end{alignedat}\right.
\end{align*}
In particular, we see that $a,b,c,\lambda$ are are different from each other and different from $0$ and $1$.

Then, we obtain if $F_{4,2}=(F_1,F_2,F_3)$ that
\tiny
\begin{align*}
	&F_1(z)=
	\frac{\splitfrac{\splitfrac{\splitfrac{\Big(44139742456659840 \, \sqrt{15} z^{8} \sqrt{215208 \, \sqrt{15} + 833497} - 56989981110573370980 \, \sqrt{15} z^{8} + 170952487440528015 \, z^{8} \sqrt{215208 \, \sqrt{15} + 833497}}{ - 220721247741925934685 \, z^{8} + 888442775916674472 \, \sqrt{15} z^{4} \sqrt{215208 \, \sqrt{15} + 833497} - 1147091808862088194764 \, \sqrt{15} z^{4} + 3440924075183568147\times  }}{\, z^{4} \sqrt{215208 \, \sqrt{15} + 833497} - 4442667472293808999113 \, z^{4} - 303994657297853240 \, \sqrt{15} \sqrt{215208 \, \sqrt{15} + 833497} + 392495488484797997540 \, \sqrt{15}}}{ - 1177366245050616585 \, \sqrt{215208 \, \sqrt{15} + 833497} + 1520128490363167610315\Big) {\left(z^{4} + 8 \, \sqrt{15} + 31\right)}}}{\splitfrac{\splitfrac{45 \, \Big(980883165703552 \, \sqrt{15} z^{4} \sqrt{215208 \, \sqrt{15} + 833497} - 1266444024679408244 \, \sqrt{15} z^{4} + 3798944165345067 \, z^{4} \sqrt{215208 \, \sqrt{15} + 833497}}{ - 4904916616487242993 \, z^{4} + 60798931459570648 \, \sqrt{15} \sqrt{215208 \, \sqrt{15} + 833497} - 78499097696959599508 \, \sqrt{15}}}{ + 235473249010123317 \, \sqrt{215208 \, \sqrt{15} + 833497} - 304025698072633522063\Big)}}\\
	&F_2(z)=\frac{\splitfrac{\splitfrac{\splitfrac{\Big(44139742456659840 \, \sqrt{15} z^{8} \sqrt{215208 \, \sqrt{15} + 833497} - 56989981110573370980 \, \sqrt{15} z^{8} + 170952487440528015 \, z^{8} \sqrt{215208 \, \sqrt{15} + 833497}}{ - 220721247741925934685 \, z^{8} + 888442775916674472 \, \sqrt{15} z^{4} \sqrt{215208 \, \sqrt{15} + 833497} - 1147091808862088194764 \, \sqrt{15} z^{4} + 3440924075183568147\times}}{ z^{4} \sqrt{215208 \, \sqrt{15} + 833497} - 4442667472293808999113 \, z^{4} - 303994657297853240 \, \sqrt{15} \sqrt{215208 \, \sqrt{15} + 833497} + 392495488484797997540 \, \sqrt{15}}}{ - 1177366245050616585 \, \sqrt{215208 \, \sqrt{15} + 833497} + 1520128490363167610315\Big) {\left(z^{4} + 8 \, \sqrt{15} + 31\right)}}}{\splitfrac{\splitfrac{45 \, \Big(980883165703552 \, \sqrt{15} z^{4} \sqrt{215208 \, \sqrt{15} + 833497} - 1266444024679408244 \, \sqrt{15} z^{4} + 3798944165345067 \, z^{4} \sqrt{215208 \, \sqrt{15} + 833497}}{ - 4904916616487242993 \, z^{4} + 60798931459570648 \, \sqrt{15} \sqrt{215208 \, \sqrt{15} + 833497} - 78499097696959599508 \, \sqrt{15}}}{ + 235473249010123317 \, \sqrt{215208 \, \sqrt{15} + 833497} - 304025698072633522063\Big)}}\\
	&F_3(z)=-\frac{\splitfrac{2 \, \Big(699302 \, \sqrt{15} z^{4} \sqrt{215208 \, \sqrt{15} + 833497} - 905413342 \, \sqrt{15} z^{4} + 2708385 \, z^{4} \sqrt{215208 \, \sqrt{15} + 833497} - 3506650795 \, z^{4}}{ - 4264614 \, \sqrt{15} \sqrt{215208 \, \sqrt{15} + 833497} + 5521560662 \, \sqrt{15} - 16516779 \, \sqrt{215208 \, \sqrt{15} + 833497} + 21384912489\Big) {\left(z^{4} + 8 \, \sqrt{15} + 31\right)} z^{2}}}{\splitfrac{3 \, \Big(699302 \, \sqrt{15} z^{4} \sqrt{215208 \, \sqrt{15} + 833497} - 905413342 \, \sqrt{15} z^{4} + 2708385 \, z^{4} \sqrt{215208 \, \sqrt{15} + 833497} - 3506650795 \, z^{4}}{
	 + 43345442 \, \sqrt{15} \sqrt{215208 \, \sqrt{15} + 833497} - 56121019962 \, \sqrt{15} + 167876175 \, \sqrt{215208 \, \sqrt{15} + 833497} - 217355775685\Big)}}
\end{align*}
\normalsize
so we see that $F_1$, $F_2$, and $F_3$ are \emph{not} polynomial functions. Indeed, we see for example that the exists two degree $4$ polynomials $P$, $Q$ such that
\begin{align*}
	F_3(z)=-\frac{2}{3}z^2(z^4+8\sqrt{15}+31)\frac{P(z)}{Q(z)}
\end{align*}
and also the coefficients of all terms in $z,z^2,z^3$ and $z^4$ of $P$ and $Q$ coincide, the constant terms of $P$ and $Q$ are respectively
\tiny
\begin{align*}
	 - 4264614 \, \sqrt{15} \sqrt{215208 \, \sqrt{15} + 833497} + 5521560662 \, \sqrt{15} - 16516779 \, \sqrt{215208 \, \sqrt{15} + 833497} + 21384912489=1.19\cdots\times 10^8
\end{align*}
\normalsize
and
\tiny
\begin{align*}
	43345442 \, \sqrt{15} \sqrt{215208 \, \sqrt{15} + 833497} - 56121019962 \, \sqrt{15} + 167876175 \, \sqrt{215208 \, \sqrt{15} + 833497} - 217355775685=-1.21\cdots\times 10^9
\end{align*}
\normalsize
which are distinct real number. Therefore $f_{4,2}:\P^1\rightarrow\C^3$ is not a null curve with simple poles at $a_1,\cdots, a_9$ and $\phi_{4,2}=\Re(f_{4,2})$ cannot be a minimal surface with $2n+1=9$ embedded planar ends.

\nocite{}
\bibliographystyle{plain}
\bibliography{biblio}

\begin{thebibliography}{10}

\bibitem{bryant}
{Robert L.} Bryant.
\newblock A duality theorem for {W}illmore surfaces.
\newblock {\em J. {D}ifferential {G}eom., 20, 23-53}, 1984.

\bibitem{bryant3}
{Robert L}. Bryant.
\newblock {\em Surfaces in conformal geometry, \emph{in The mathematical
  heritage of Hermann Weyl}}.
\newblock Proceedings of {S}ymposia in {P}ure {M}athematics, 1987.

\bibitem{bryantcurves}
Robert~L. Bryant.
\newblock Projective, contact, and null curves.
\newblock {\em Preprint, 16th December}, 2017.

\bibitem{calabi}
Eugenio Calabi.
\newblock Minimal immersions of surfaces in {E}uclidean spheres.
\newblock {\em J. of Differential Geom., 1, 111-125}, 1967.

\bibitem{costa2}
{Celso J.} Costa.
\newblock Complete minimal surfaces in $\mathbb{R}^3$ of genus one and four
  planar embedded ends.
\newblock {\em Proc. Amer. Math. Soc. 119, no. 4, 1279–1287}, 1993.

\bibitem{dierkes}
Ulrich Dierkes, Stefan Hildebrandt, and Friedrich Sauvigny.
\newblock {\em Minimal surfaces}.
\newblock SPringer-Verlag, Grundlehren der mathematischen Wissenschaften, vol.
  339, 2010.

\bibitem{griffiths}
Philipp Griffiths and Joseph Harris.
\newblock {\em Principles of Algebraic Geometry}.
\newblock Pure and Applied Mathematics. Wiley-Interscience, New York, 1978.

\bibitem{helicoid}
David Hoffman, Martin Traizet, and Brian White.
\newblock Helicoidal minimal surfaces of prescribed genus.
\newblock {\em Acta Math. 216, no. 2, 217–323.}, 2016.

\bibitem{jorge}
{Luquesio P.} Jorge and {William H.}~Meeks III.
\newblock The topology of complete minimal surfaces of finite total {G}aussian
  curvature.
\newblock {\em Topology, Vol. 22n Ni. 2, pp. 203-221}, 1983.

\bibitem{klein}
Felix Klein.
\newblock \"{U}ber die {T}ransformation der allgemeinen {G}leichung des zweiten
  {G}rades zwischen {L}inien-{C}oordinaten auf eine canonische {F}orm.
\newblock {\em Mathematische Annalen 23: 539.
  https://doi.org/10.1007/BF01446603}, 1884.

\bibitem{kollar}
János Kollár.
\newblock {\em Lectures on resolution of singularities}.
\newblock Annals of Mathematics Studies, 166. Princeton University Press,
  Princeton, NJ. vi+208 pp. ISBN: 978-0-691-12923-5; 0-691-12923-1, 2007.

\bibitem{kusnerpacific}
Robert Kusner.
\newblock Comparison surfaces for the {W}illmore problem.
\newblock {\em Pacific {J}. {M}ath., {V}ol. 138, No. 2}, 1989.

\bibitem{linriv}
Fang-Hua Lin and Tristan Rivière.
\newblock Energy quantization for harmonic maps.
\newblock {\em Duke Math. J. Volume 111, Number 1, 177-193}, 2002.

\bibitem{classification}
Alexis Michelat and Tristan Rivi{\`e}re.
\newblock The {C}lassification of {B}ranched {W}illmore {S}pheres in the
  $3$-{S}phere and the $4$-{S}phere.
\newblock {\em arXiv:1706.01405}, 2017.

\bibitem{sagepaper}
Alexis Michelat and Tristan Rivière.
\newblock Computer-{A}ssisted {P}roof of the {M}ain {T}heorem of\\ '{T}he
  {C}lassification of {B}ranched {W}illmore {S}pheres in the $3$-{S}phere and
  the $4$-{S}phere'.
\newblock {\em arXiv:1711.10441}, 2017.

\bibitem{blow-up}
Alexis Michelat and Tristan Rivière.
\newblock Higher {R}egularity of {W}eak {L}imits of {W}illmore {I}mmersions
  {I}.
\newblock {\em arXiv:1904.04816}, 2019.

\bibitem{blow-up2}
Alexis Michelat and Tristan Rivière.
\newblock Higher {R}egularity of {W}eak {L}imits of {W}illmore {I}mmersions
  {II}.
\newblock {\em arXiv:1904.09957}, 2019.

\bibitem{osserman}
Robert Osserman.
\newblock On {C}omplete {M}inimal {S}urfaces.
\newblock {\em Arch. Rational Mech. Anal., 13, p. 392–404}, 1963.

\bibitem{peng0}
Chia-Kuei Peng.
\newblock Some new examples of minimal surfaces in $\mathbb{R}^3$.
\newblock {\em Preprint, MSRI}, 1986.

\bibitem{peng1}
Chia-Kuei Peng and Liang Xiao.
\newblock Willmore {S}urfaces and {M}inimal {S}urfaces with {F}lat {E}nds.
\newblock {\em Geometry and Topology of {S}ubmanifolds}, 2000.

\bibitem{schoenPlanar}
Richard~M. Schoen.
\newblock {U}niqueness, {S}ymmetry, and {E}mbeddedness of {M}inimal {S}urfaces.
\newblock {\em J. {D}ifferential {G}eom., 18, 791-809}, 1983.

\bibitem{shamaev}
{Èlly I.} Shamaev.
\newblock On a family of minimal tori in $\mathbb{R}^3$ with planar ends.
\newblock {\em Sibirsk. Mat. Zh. 46, no. 6, 1407--1426; translation in Siberian
  Math. J. 46, no. 6, 1135–1152}, 2005.

\bibitem{peng2}
Liang Xiao.
\newblock On complete minimal surfaces with parallel and flat ends.
\newblock {\em In: Jiang B., Peng CK., Hou Z. (eds) Differential Geometry and
  Topology. Lecture Notes in Mathematics, vol 1369. Springer, Berlin,
  Heidelberg}, 1989.

\end{thebibliography}
\end{document}